\documentclass{article}

\usepackage{cite}
\usepackage{graphicx} % Required for inserting images
\usepackage{booktabs}

\usepackage{amsmath,amsfonts,amssymb,amsthm}
\newcommand\real{\mathbb{R}}

\DeclareMathOperator{\sat}{sat}
\DeclareMathOperator{\dz}{dz}
\DeclareMathOperator{\He}{He}

\newenvironment{smatrix}%          environment name
{\left[\begin{smallmatrix}}%            begin code
{\end{smallmatrix}\right]}% 

\newcommand\smallmat[1]{\left[\begin{smallmatrix}#1\end{smallmatrix}\right]}
\newcommand\bigmat[1]{\begin{bmatrix}#1\end{bmatrix}}

\newcommand\KT{K^\top}

\newtheorem{prop}{Proposition}
\newtheorem{coro}{Corollary}

\newtheorem{theo}{Theorem}
\newtheorem{lemma}{Lemma}

\newtheorem{remm}{Remark}
\newtheorem{exx}{Example}

\newenvironment{rem}{\begin{remm}\rm }{\hfill \hspace*{1pt} \hfill $\circ$\end{remm}}
\newenvironment{example}{\begin{exx}\rm }{\hfill \hspace*{1pt} \hfill $\lrcorner$\end{exx}}

\def\BibTeX{{\rm B\kern-.05em{\sc i\kern-.025em b}\kern-.08em
    T\kern-.1667em\lower.7ex\hbox{E}\kern-.125emX}}

\newcommand\Acl{A_{\mathrm{c\ell}}}

\begin{document}
% \title{Achieving GAS with non-exponentially stabilizable saturated linear plants}
% \title{Globally stabilizing saturated linear feedback for non-exponentially stabilizable linear plants}
% \title{Nonquadratic global asymptotic stability with saturated linear feedback of asymptotically null-controllable plants}
\title{Nonquadratic global asymptotic stability certificates for saturated linear feedbacks}

\author{Andrea Cristofaro\thanks{{A. Cristofaro is with the Department of Computer, Control and Management Engineering, Sapienza University of Rome, Italy (e-mail: andrea.cristofaro@uniroma1.it  )}}\, and Luca Zaccarian\thanks{L. Zaccarian is with LAAS-CNRS, Toulouse, France and the University of Trento, Italy (e-mail: zaccarian@laas.fr)}
}

\date{}
\maketitle

\begin{abstract}
We establish sufficient conditions for positive (semi-)definiteness, with or without radial unboundedness, for nonquadratic Lyapunov function constructed as sign-indefinite quadratic forms involving the state and the deadzone of a suitable input. We then use these conditions to build weak nonquadratic Lyapunov functions establishing global asymptotic stability of linear systems in feedback through a saturation, leveraging invariance principles. Our results
are shown to be non-conservative (necessary and sufficient)
for a family of well known prototypical examples of linear SISO feedbacks that are not globally exponentially stabilizable (the so-called ANCBI plants). Our multi-input extension leads to convex stability analysis tests, formulated as linear matrix inequalities that are applicable to ANCBI 
non-globally-exponentially-stabilizable plants. 
\end{abstract}

\section{Introduction}

It has been established in the 1980s \cite{Sontag84}, \cite{Schmitendorf80}, \cite{Lasserre92} that a linear plant cannot be globally exponentially or asymptotically stabilized by bounded inputs when it contains
exponentially unstable dynamics. 
More precisely, with continuous-time dynamics, those works prove that 
a linear plant subject to a decentralized, possibly non-symmetric, input saturation function
\begin{align}
\label{eq:plant0}
\dot{x}=A_0 x-B\sat(u),
\end{align}
(where $x\in\mathbb{R}^n,\ u\in\mathbb{R}^m$ and $\sat(u)=u$ in a small enough neighborhood of $u=0$)
can be global asymptotically stabilized (GAS) by a feedback
$u=\gamma(x)$ if and only if $A_0$ only has eigenvalues with non-positive real part. Moreover, global exponential stability (GES) can be obtained if and only if $A_0$ is Hurwitz, namely
\begin{align}
\label{eq:HurwitzA0}
\exists P_0\succ0 :\; A_0^\top P_0 +P_0 A_0 \prec0,
\end{align}
an intrinsic limitation that also holds in the infinite-dimensional case \cite[Prop. 2]{astolfiCDC22}.

An interesting special case corresponds to plants that are not globally exponentially stabilizable (namely, $A_0$ is not Hurwitz) but can be globally asymptotically stabilized (namely, $A_0$ does not have exponentially unstable eigenvalues). This class of linear plants has been named Asymptotically Null-Controllable with Bounded Inputs (ANCBI) in the early literature and their stabilization was known to be an authentic nonlinear problem because Fuller had already established in 1969 that no saturated linear controller could globally asymptotically stabilize the triple integrator \cite{Fuller69}. In fact, early works like the
nested saturations \cite{Teel92,Sussmann94} and scheduled Riccati \cite{Megretski96,HouSaberiAuto98,saberi2000simultaneous}, managed to provide global asymptotic stabilizers for any ANCBI plant. Follow-up works recognized that the global asymptotic stabilizers from the 1990s could be improved in terms of transient speed and proposed alternative scheduled or switched nonlinear solutions leading to improved performance (see, e.g., \cite{Forni10, Zhou10} or the more recent works
\cite{ValmorbidaTAC17,LA14}, which propose polynomial and rational globally asymptotically stabilizing feedback).
Another strategy, originally formulated for the regulation problem~\cite{CLPV03,LL15}, parametrizes state feedback gain by way of a scalar, bounded nonlinear function. 
A further specific line of research addresses performance-oriented 
feedback for specific plants like the saturated double integrator \cite{Forni10, Zhou10,tyan1999global},
whose relevance stems from its ubiquitous presence in industrial systems, as emphasized in the survey work \cite{rao2001naive}.

While nonlinear feedbacks remain a valid approach, there are still many relevant cases, such as those of non-polynomially unstable plants, where linear saturated feedbacks are known to perform well. In this case, the closed loop between plant \eqref{eq:plant0} and a linear feedback $u=Kx$
with $K\in\mathbb{R}^{m\times n}$,
corresponds to 
\begin{align}
\label{eq:linearCL}
\dot{x}=A_0 x-B\mathrm{sat}(Kx) = (A_0-BK)x+B\dz(Kx),
\end{align}
where we introduced the deadzone function $\dz(u) := u-\sat(u)$, taking its name from the fact that it is zero in a neighborhood of the origin. In fact, close to zero, the closed loop \eqref{eq:linearCL} is linear, which means that asymptotic stability of the origin holds only if $\Acl:=A_0-BK$ is Hurwitz. Equation \eqref{eq:linearCL} represents a quite general setting, including dynamic output feedback controllers with strictly proper linear plants (see, e.g., \cite[Remark 2]{ValmorbidaTAC22}).

When $A_0$ is not Hurwitz, typical quadratic Lyapunov approaches are not viable for establishing global asymptotic stability of the origin for \eqref{eq:linearCL}, because the ensuing quadratic conditions would imply global exponential stability (which is not achievable from a bounded input). A simple example illustrating this fact is the scalar system $\dot x = -\sat(k x)$, $k>0$ for which a Lyapunov function $V(x)=kx^2$ satisfies $\dot V(x) =-kx\sat(kx)$, which is clearly negative definite, but whose negativity cannot be established by quadratic bounds. 
A promising Lyapunov structure for capturing the key features of the simple scalar example corresponds to the extended quadratic forms
$\begin{smatrix}
    x\\
    \dz(Kx)
\end{smatrix}^\top
\smallmat{Q_{11} & Q_{12}\\ Q_{12}^\top & Q_{22}} 
\begin{smatrix}   x\\   \dz(Kx)\end{smatrix}$
first proposed in  \cite{pri:gia/tac2001,dai:hu:tee:zac/scl2009}, insisting on positive definiteness of matrix $Q=
\smallmat{Q_{11} & Q_{12}\\
Q_{12}^\top & Q_{22}}$. Later works   \cite{val:gar:zac/auto2017,val:dru:dun/tac2019,dru:val:dun/csl2018,LiLinIFAC2015,ValmorbidaTAC22}
recognized that positive definiteness of $Q$ is not needed for positive definiteness of the quadratic form, and led to less conservative stability analysis conditions, and the follow-up designs
in \cite{PriuliLCSS22,CalderonLCSS24} and references therein. Nevertheless, all of the above-cited works insist on (local or global) exponential stability properties by requiring the same homogeneous bounds (typically quadratic) for the Lyapunov function and its directional derivative. 
A different and more sophisticated approach is required in the recent work 
\cite{ChitourMCSS23}, which
 constructs a converging family of stabilizers
to prove that the saturated double oscillator can be stabilized by a specific linear saturated feedback. The ensuing construction does not violate the necessary condition in \cite{Fuller69} (see also \cite{sussmannCDC91}).

In this paper, we investigate relaxed positivity and radial unboundedness conditions for these so-called {\em sign-indefinite quadratic forms} allowing us to draw conclusions about global asymptotic stability
of the origin for \eqref{eq:linearCL}
even for cases where global exponential stability cannot be obtained, due to the above mentioned limitations highlighted in \cite[Prop. 2]{astolfiCDC22}.
Our results hinge upon characterizing conditions on $Q$ ensuring positivity without imposing quadratic bounds.
Our construction is effective for matrices $A_0$ satisfying
the following relaxed version of \eqref{eq:HurwitzA0}:
\begin{align}
\label{eq:LiuAkasaka}
\exists P_0\succeq0 :\; -S_0:=A_0^\top P_0 +P_0 A_0 \preceq 0,\;  \ker P_0 \subseteq \ker A_0,
\end{align}
which resembles the condition proposed in \cite[eqn. (3)]{LA14} with a subtle yet important difference.~\footnote{In particular, in \cite[eqn. (3)]{LA14} the inclusion between $\ker P_0$ and $\ker A_0$ is reversed as compared to \eqref{eq:LiuAkasaka}. Reversing the inclusion appears to be inevitable, otherwise, with \cite[eqn. (3)]{LA14}, the trivial selection $P_0=0$ would always satisfy the constraint.} In fact, we show in the appendix that the following result holds, which clarifies the class of systems that we can address with our construction.
This result was suggested (without any proof, and referring to a technical report that is no longer available) in \cite[\S 2]{LA14}. The subtle difference between \eqref{eq:LiuAkasaka} and \cite[eqn. (3)]{LA14} motivated us to provide a full proof in the appendix.

\begin{prop} \label{prop:LiuAkasaka}
Condition \eqref{eq:LiuAkasaka} holds if and only if $A_0$ has eigenvalues in the closed left-half plane with Jordan blocks of dimension at most 2 for the eigenvalues at the origin and  simple Jordan blocks for the eigenvalues in the rest of the imaginary axis.
\end{prop}

% While Proposition~\ref{prop:LiuAkasaka} was suggested (without any proof, and referring to a technical report that is no longer available) in Section 2 of \cite{LA14}, the subtle difference between \eqref{eq:LiuAkasaka} and \cite[eqn. (3)]{LA14} motivated us to provide a full proof, which is reported in the Appendix.

It is worth commenting on the fact that, due to Proposition~\ref{prop:LiuAkasaka}, the saturated double oscillator described in \cite{ChitourMCSS23} does not satisfy \eqref{eq:LiuAkasaka}
and is therefore not a system that we manage to stabilize with the linear state feedback characterizations in this paper. A relevant case exhibiting such a saturated double oscillator structure corresponds to the (well known to few) example proposed by E. Sontag in the 1990s, corresponding to the feedback \eqref{eq:linearCL} with 
$A_0 = \smallmat{0& 1& 0& 0\\ 0& 0 &1& 0\\ 0 &0& 0& 1\\ -1& 0& -2& 0}$, $B^\top = \bigmat{0& 0 & 0& 1}$, $K = \bigmat{1&1&2&1}$. While we acknowledge that our construction does not address this challenge, we believe that our tools provide an important step towards understanding the intricate behavior of these saturated linear feedbacks.

% . Matrices satisfying \eqref{eq:LiuAkasaka} include single and double integrators, and matrices having single eigenvalues on the imaginary axis.

The paper is organized as follows.
Section~\ref{sec:nonquadratic_lemmas} provides novel results on sign-indefinite quadratic forms that are positive definite, radially unbounded, but not quadratically bounded.
Section~\ref{sec:GASresult} applies those results to relaxed GAS conditions that do not imply GES.
Section~\ref{sec:SI} uses such conditions on
low-dimensional ANCBI prototypical examples where we build universal Lyapunov certificates.
Then, the construction is generalized to higher order plants satisfying \eqref{eq:LiuAkasaka} in
Section~\ref{sec:generalization}, where
we
propose analysis and design conditions based on linear matrix inequalities, ensuring global asymptotic stability, even though global exponential stability is not achievable.
 Numerical examples are discussed in 
Section~\ref{sec:examples} and then conclusions are drawn in 
Section~\ref{sec:conclusions}.

% An insightful characterization of properties of $A_0$ that would allow for this type of bounds,
% while accounting for the intrinsic limitations of saturated linear stabilizers highlighted by Fuller in \cite{Fuller69},
% has been proposed in \cite{LA14}, and corresponds to 

% An interesting open question is whether \eqref{eq:LiuAkasaka} together with some additional condition on the input matrix $B$ is enough for ensuring 

{\bf Notation}. $\mathbb{D}_{\geq0}^{m}$ (respectively, $\mathbb{D}_{>0}^{m}$) is the set of positive semi-definite (respectively, positive definite) diagonal matrices of dimension $m\times m$. Similarly, 
$\mathbb{S}_{\geq0}^{n}$ (respectively, $\mathbb{S}_{>0}^{n}$) is the set of positive semi-definite (respectively, positive definite) symmetric matrices of dimension $n\times n$.
Given a square matrix $Q$, $\sigma(Q)$ denotes the set of its eigenvalues and $\He(Q) = Q+Q^\top$.

\newcommand\myQ{Q}

\section{Nonquadratic positivity of sign-indefinite extended forms}
\label{sec:nonquadratic_lemmas}

In this section we specify sufficient conditions for positive definiteness and possibly radial unboundedness of the sign-indefinite quadratic forms 
\begin{equation}\label{eq:quadratic_form}
    Y(x)=\begin{bmatrix}
    x\\
    \dz(Kx)
\end{bmatrix}^\top \underbrace{\begin{bmatrix}
Q_{11} & Q_{12}\\
Q_{12}^\top & Q_{22}
\end{bmatrix}}_{=\myQ }\begin{bmatrix}
    x\\
    \dz(Kx)
\end{bmatrix}
\end{equation}
proposed in \cite{pri:gia/tac2001,dai:hu:tee:zac/scl2009} and characterized in \cite{ValmorbidaTAC22} when insisting on sign-indefinite matrices $Q$. While  \cite{ValmorbidaTAC22} provided conditions for establishing quadratic upper and lower bounds on function $Y$, we relax those conditions here to be able to exploit $Y$ for establishing asymptotic, rather than exponential, stability properties. 
% We first state our main result with nonquadratic bounds, and then we present several illustrative examples.

%
% Let us consider a linear system controlled by a saturated feedback
% \begin{align}
% \label{eq:plant}
% \dot{x}=A_0 x-B\mathrm{sat}(Kx)
% \end{align}
% where $x\in\mathbb{R}^n,\ u\in\mathbb{R}^m$ and with $K\in\mathbb{R}^{m\times n}$ such that
% $$
% \sigma(A_0 -BK)\subset\mathbb{C}^{-}
% $$
% Following the ideas originally introduced in \cite{}, we look for a quadratic Lyapunov function of the form
%
% where $dz:=Kx-\mathrm{sat}(Kx)$, and with $\myQ $ possibly being sign-indefinite. Let us also 
Let us first summarize a few properties.
 % known as sector conditions. 
 Given an arbitrary non-negative diagonal matrix $T_0\in\mathbb{D}_{\geq0}^{m}$ one has the sector condition
$$
2\dz(u)^\top T_0 \sat(u)\geq 0.
$$
When $u=Kx$, this sector condition reads as the quadratic form
 % can be written, respectively, as
\begin{align}
\label{eq:lower_cond}
\begin{bmatrix}
    x\\
    \dz(Kx)
\end{bmatrix}^\top \underbrace{\begin{bmatrix}
0&K^\top T_0\\
T_0 K&-2T_0
\end{bmatrix}}_{=\Sigma_0}\begin{bmatrix}
    x\\
    \dz(Kx)
\end{bmatrix}&\geq0,
% \\
% \label{eq:upper_cond}
% \begin{bmatrix}
%     x\\
%     \dz(Kx)
% \end{bmatrix}^\top\underbrace{\begin{bmatrix}
% 2K^\top T_1 K &-K^\top T_1\\
% -T_1 K&0
% \end{bmatrix}}_{=\Sigma_1}\begin{bmatrix}
%     x\\
%     \dz(Kx)
% \end{bmatrix}&\geq0
\end{align}
% We will refer to the conditions \eqref{eq:lower_cond} and \eqref{eq:upper_cond}, respectively, as {\it lower} and {\it upper} sector condition (due to the lower or upper position of the nontrivial diagonal term).
Furthermore, given an arbitrary positive semi-definite matrix $R\in\mathbb{R}^{m\times m}$ 
define the following matrix, which is obviously positive semi-definite,
\begin{align}
\label{eq:R_cond}
\Sigma_{R}:=\begin{bmatrix}
K^\top R K & -K^\top R\\
-R K & R
\end{bmatrix}\succeq 0.
\end{align}
The following theorem is our first main result.

\begin{theo}\label{theo:definite}
Consider function $x \mapsto Y(x)$ in \eqref{eq:quadratic_form} associated with matrix $Q=Q^\top$, and the matrices $\Sigma_0$,  $\Sigma_R$ in \eqref{eq:lower_cond},   \eqref{eq:R_cond}. Then the following sufficient conditions hold true.
\begin{enumerate}
\item[(i)] Assume that $T_0\in\mathbb{D}_{\geq0}^{m}$ and $R \in\mathbb{S}_{\geq 0}^{m}$ exist such that
\begin{align}
\label{eq:Th1_psd}
Q -\Sigma_0-\Sigma_R\succeq 0
\end{align}
Then function $Y$ is positive semi-definite.

\item[(ii)] Assume that $T_0\in\mathbb{D}_{\geq0}^{m}$ and $R\in\mathbb{S}_{>0}^{m}$ exist such that
\begin{align}
\label{eq:Th1_pd_LowerBound}
Q  -\Sigma_0- \Sigma_R\succeq 0 ,
\end{align}
then $Y$ is positive semi-definite and enjoys the global lower bound $Y(x) \geq \lambda_{\min}(R)|\mathrm{sat}(Kx)|^2$.

\item[(iii)] Assume that $R\in\mathbb{S}_{>0}^{m}$ exists such that
\begin{align}
\label{eq:Th1_pd}
Q_{11}\succ 0 \mbox{ and } Q -\Sigma_R\succeq 0 ,
\end{align}
then $Y$ is positive definite, but not necessarily radially unbounded.

\item[(iv)] Assume that $T_0\in\mathbb{D}_{>0}^{m}$ exists such that
\begin{align}
\label{eq:Th1_pdru}
Q_{11}\succ 0 \mbox{ and } Q -\Sigma_0\succeq 0
\end{align}
then $Y$ is positive definite and radially unbounded.
\end{enumerate}
\end{theo}

\begin{proof}
To prove item (i), exploiting \eqref{eq:lower_cond} and \eqref{eq:R_cond}, we have for any $x\in \real^n$,
\begin{align*}
    Y(x) &\geq Y(x) - \begin{smatrix}
    x\\
    \dz(Kx)
\end{smatrix}^\top(\Sigma_0+\Sigma_R)\begin{smatrix}
    x\\
    \dz(Kx)
\end{smatrix}\\
  &= \begin{smatrix}
    x\\
    \dz(Kx)
\end{smatrix}^\top(Q-\Sigma_0-\Sigma_R)\begin{smatrix}
    x\\
    \dz(Kx)
\end{smatrix}\geq 0,
\end{align*}
where the last inequality follows from \eqref{eq:Th1_psd}.

To prove item (ii), it is straightforward to check the equality 
$$
\begin{bmatrix}
x\\
\dz(Kx)
\end{bmatrix}^\top\Sigma_{R}\begin{bmatrix}
x\\
\dz(Kx)
\end{bmatrix}=\sat(Kx)^\top R\, \sat(Kx).
$$
Then, using \eqref{eq:lower_cond}, condition $\myQ -\Sigma_0 -\Sigma_R \succeq 0$ in \eqref{eq:Th1_pd_LowerBound} implies
\begin{align}
Y(x)\geq\sat(Kx)^\top R\, \sat(Kx)\geq \lambda_{\min}(R)|\mathrm{sat}(Kx)|^2,
\label{eq:inProofThm1}
\end{align}
which proves item~(ii).

Let us continue with item (iii)
and note that  \eqref{eq:inProofThm1} holds because of the assumption $Q-\Sigma_R\succeq 0$ in \eqref{eq:Th1_pd}.
Then,
we can observe that the right-hand side of \eqref{eq:inProofThm1} is positive for any $x$ such that $Kx\neq 0$. On the other hand, for any $x\neq 0$ such that $Kx= 0$, we have $\dz(Kx)=0$ and, due to the assumption $Q_{11}\succ 0$ in \eqref{eq:Th1_pd},
$$
Y(x)= 
\begin{bmatrix}
x\\
\dz(Kx)
\end{bmatrix}^\top \myQ \begin{bmatrix}
x\\
\dz(Kx)
\end{bmatrix}=x^\top Q_{11} x>0,
$$
thus proving the positive definiteness of $Y$.

Item (iv) is proven by establishing a contradiction. 
First of all, note that item (i) holds, therefore $Y$ is positive semi-definite. To show positive definiteness,
suppose for contradiction
that $x_0 \neq 0$
exists such that $Y(x_0)=0$, and denote
$\xi_0=[x_0^\top \; \dz(Kx_0)^\top]^\top\neq 0$. 
Since $Q_{11}\succ 0$ by assumption, the previous condition may only be satisfied for $\dz(Kx_0)\neq0$. On the other hand $T_0\succ 0$, 
therefore $\dz(Kx_0) \neq 0$ implies $2\dz(Kx_0)^\top T_0\ \mathrm{sat}(Kx_0) = \xi_0^\top \Sigma_0 \xi_0 > 0$ (where we used \eqref{eq:lower_cond}). Combining this last bound with the second bound in \eqref{eq:Th1_pdru},
we obtain
$$
Y(x_0)=\xi_0^\top\myQ \xi_0\geq \xi_0^\top\Sigma_0\xi_0>0,
$$
thus proving that $Y$ is actually positive definite. Let us now prove radial unboundedness. Again, for contradiction, assume that $\tilde{x}\neq 0$ and $M>0$ exist such that 
\begin{align}
\label{eq:contraTh1item2}
\lim_{\alpha\rightarrow +\infty} Y(\alpha \tilde{x})\leq M
\end{align}
Since $Q_{11}\succ0$, then for any $x$ satisfying $Kx=0$ the quadratic form is trivially radially unbounded. Therefore, it is necessary that $dz(K\tilde{x})\neq0$ for \eqref{eq:contraTh1item2} to hold. For simplicity suppose also that $K_{i} \tilde{x}>0$ for some $i\in\{1,...,m\}$ (the case 
$K_{i} \tilde{x}<0$ follows parallel steps). Evaluating the quadratic form at $\alpha\tilde{x}$ with $\alpha>0$,
exploiting \eqref{eq:Th1_pdru}
and using $T_0\succ 0$ yields
\begin{align}
\label{eq:Thm1_last_step}
Y(\alpha\tilde{x})&\geq \begin{bmatrix}
\alpha \tilde{x}\\
\dz(K\alpha\tilde{ x})\end{bmatrix}^\top \Sigma_0 \begin{bmatrix}
\alpha \tilde{x}\\
\dz(K\alpha\tilde{ x})\end{bmatrix}\\
\nonumber
&\geq 2\dz(K_i\alpha\tilde{x})T_{0_i} \sat(K_i\alpha\tilde{x})
\end{align}
For $\alpha$ sufficiently large so that $K_i\alpha\tilde{x}$ exceeds the saturation limit, the expression at the right-hand side of \eqref{eq:Thm1_last_step} becomes
$2\alpha T_{0_i}K_i\tilde{x}-2T_{0_i}$
which, in turn, implies that, 
for any
$\alpha>\dfrac{M+2T_{0_i}}{2T_{0_i}K_i\tilde{x}}$,
we have $Y(\alpha\tilde{x})>M$, thus contradicting \eqref{eq:contraTh1item2}.
\end{proof}

We emphasize that the results in Theorem~\ref{theo:definite}
characterize a rather specific selection of the multipliers and of the matrices characterizing the extended quadratic form \eqref{eq:quadratic_form}. We clarify below, by some examples, that a nontrivial richness of behavior can be obtained by suitable selections.

\begin{example}
Item (iii) truly allows for functions that are positive definite by not radially unbounded. For example, 
given any full row rank gain $K$,
consider
selecting $Q_{12}=-K^\top R, Q_{22}=R$ and $Q_{11}=K^\top R K+X_{11}$ with 
$$X_{11}=(I-K^\top(KK^\top)^{-1}K)^\top W(I-K^\top(KK^\top)^{-1}K)$$
for any $W=W^\top\succ 0$. Then,
observing that matrix $\Pi_K := (I-K^\top(KK^\top)^{-1}K)$, appearing in the above selection of $X_{11}$, is a projection matrix (satisfying $\Pi_K K^\top=0$ and $\Pi_K K^\perp =  K^\perp$, where $K^\perp$ denotes the orthogonal complement of $K^\top$) we may conclude that
\eqref{eq:Th1_pd} holds because $KQ_{11}K^\top=KK^\top RKK^\top \neq 0$ and $(K^\perp)^\top Q_{11}K^\perp= (K^\perp)^\top W K^\perp \neq 0$, so that $Q_{11}$ is positive definite. Moreover, $Q-\Sigma_R = \smallmat{Q_{11} & 0\\0 &0}\succ 0$.

However, picking $x= K^\top y$ for any nonzero $y\in\mathbb{R}^m$, one has $X_{11}x=0$ and therefore, for any $\lambda \in \real$,
$$
Y(\lambda x)=\sat(\lambda KK^\top y)^\top R\, \sat(\lambda KK^\top y),
$$
which is uniformly bounded for any arbitrarily large value of $\lambda$, thereby showing that $Y$ is not radially unbounded. 
% $$
% \lim_{\lambda\rightarrow +\infty} V(\lambda \tilde{x})\leq \lambda_{\max}(R)\, m.
% $$
\end{example}

\begin{example}\label{example:SPD1}
 It may seem intuitive that positive definiteness in item (iv) of Theorem~\ref{theo:definite} could also follow for $T_0\in\mathbb{D}_{\geq0}^{m}$.
Let us show, with a counter-example, that 
positive definiteness of $T_0$ is generally necessary.
% the condition $\myQ -
% \Sigma_0\succeq 0$, with $T_0\in\mathbb{D}_{\geq0}^{m}$ and $P_{11}\succ 0$, is not sufficient for positive definiteness. 
Consider the selection
$$
\myQ =\begin{smatrix}
1 & 2 &0\\
2 & 5  &-1\\
0&-1&1
\end{smatrix},\quad P_{11}=\begin{smatrix}
1 & 2\\
2 & 5 
\end{smatrix}\succ 0,
$$
and observe that $\sigma(\myQ )=\{6,1,0\}$, thus implying that condition \eqref{eq:Th1_pdru} holds with $T_0=0$. 

Consider now the feedback matrix $K=[1\ \tfrac{7}{2}]$ and the point $x^\star=[-4\ 2]^\top$. It is easy to compute
$\dz(Kx^\star)=Kx^\star-\sat(K x^\star)=(-4+7)-1=2$,
so that 
$\xi^\star=[{x^\star}^\top \dz(Kx^\star)]^\top = [-4\ 2\ 2]^\top$.
% $$
% dz^\star=K_1x^\star-\mathrm{sat}(K_1x^\star)=(-4+7)-1=2
% $$
% and
% $$
% \xi^\star=[{x^\star}^\top dz^\star]^\top
% $$
It is immediate to check that $\myQ \xi^\star=0$, thus showing that $Y$ is not positive definite.
\end{example}

\begin{example}
The previous example correctly suggests that in the single input case,
when 
\eqref{eq:Th1_pdru} holds with $T_0=0$ and $\myQ \succeq0$
such that $0\in\sigma(\myQ )$,
there exists a whole set of feedback gains $K$
such that
the quadratic form $Y$ in \eqref{eq:quadratic_form} is not positive definite.

This is indeed true and such gains $K$ can be constructed by 
selecting $z=\smallmat{z_{1}\\ z_2}\in\ker{\myQ }$, so that
$$
\begin{array}{l}
Q_{11}z_1+Q_{12}z_2=0\smallskip\\ Q_{12}^\top z_1+Q_{22}z_2=0
\end{array}
$$
Then we may pick any gain $K_1$ (e.g., $K_1=z_1^\top$) such that $\dz(K_1 z_1)\geq0$,
%say $dz>0$ for the sake of simplicity. When 
and $K_1 z_1>z_2>0$, so that the equation 
$\lambda z_2=\lambda K_1 z_1-1$
admits the positive solution $\lambda^*=1/(K_1 z_1-z_2)>0$. 
Then we obtain 
% This proves indeed that, for $Kz_1>z_2$, one has
$
% \myQ 
% \begin{bmatrix}
% \lambda^* z_1\\
% \lambda^* z_2
% \end{bmatrix}=
\myQ 
\begin{bmatrix}
\lambda^* z_1\\
\dz(\lambda^* K z_1)
\end{bmatrix}=0
$, which proves that $Y$ in  \eqref{eq:quadratic_form} is not positive definite.
\end{example}

\begin{example}
The previous example may further suggest that perhaps 
the conclusions of item (iv) of Theorem~\ref{theo:definite} could also be true for $T_0\in\mathbb{D}_{\geq0}^{m}$ whenever $0\notin\sigma(Q)$. However, this is not the case as illustrated by the selection below.

Consider the selection $Q_{11} = 4$, 
$Q_{12} = \smallmat{-2 & -\tfrac{1}{2}}$ 
and $Q_{22} = \smallmat{1 & 0 \\ 0 & -1}$. It is easy to check that $\det(Q) = -\tfrac{1}{4} \neq 0$, therefore $0\notin\sigma(\myQ _{ext})$. 
Select now $K = \smallmat{3\\ -0.5}$ and the positive semi-definite $T_0 = \smallmat{0 & 0\\ 0& 1}\succ 0$, so that
$\Sigma_0=\begin{bmatrix}
0 & K^\top T_0\\
T_0K & -2T_0
\end{bmatrix} = \smallmat{0 & 0 &-0.5 \\ 0 & 0 &0 \\ -0.5 &0 &-2}$.

With these selections we can check that \eqref{eq:Th1_pdru} holds because
$Q-\Sigma_0 = \smallmat{
4      &     -2    & 0\\
-2       &     1   &         0\\
0    &        0      &     1}\succeq 0$. Nonetheless, picking $x^\star=1$,
we have $\dz(Kx^\star)=\dz(\smallmat{3\\ -0.5}) = \smallmat{2\\0}$. Then it is immediate
to check that 
$\bigmat{x^\star \\ \dz(Kx^\star)}^\top Q\bigmat{x^\star \\ \dz(Kx^\star)}=0$,
thus showing that $Y$ is not positive definite.
\end{example}

\section{Lyapunov construction for GAS but not GES}
\label{sec:GASresult}

We construct in this section a Lyapunov certificate of GAS for 
\eqref{eq:linearCL} under assumptions that are mild enough to cover cases where it impossible for GES to hold for any feedback gain $K$. To this end we leverage the results in 
Theorem~\ref{theo:definite}. 

Based on \eqref{eq:quadratic_form}, consider the following block diagonal Lyapunov function candidate
\begin{equation}\label{eq:Lyapunov_f}
    V(x)=\frac{1}{2}\begin{bmatrix}
    x\\
    \dz(Kx)
\end{bmatrix}^\top \underbrace{\begin{bmatrix}
P_{11} & 0 \\
0 & P_{22}
\end{bmatrix}}_{=P}\begin{bmatrix}
    x\\
    \dz(Kx)
\end{bmatrix}
\end{equation}
where $P_{22}$ is diagonal and possibly sign indefinite.

% \red{Secondo me qui dobbiamo assumere che $P_{22}$ e' diagonale!}\\
The
block-diagonal structure 
of \eqref{eq:Lyapunov_f}
allows obtaining a directional derivative along \eqref{eq:linearCL} still enjoying the quadratic structure in \eqref{eq:quadratic_form}. This is because of the suggestive identity $\frac{d}{dt}\dz(Kx)^2 = 2\dz(Kx)^\top K\dot x$, which, together with the fact that $P_{22}$ is diagonal, allows defining, for almost all $x\in \real^n$, 
\begin{align}
\nonumber
    \dot V(x)&=x^\top P_{11}(\Acl x+B\dz) + \dz^\top P_{22}K(\Acl x+B\dz)\\
    &= \begin{bmatrix}
    x\\
    \dz
\end{bmatrix}^\top \begin{bmatrix}
  P_{11}\Acl & P_{11}B \\
 P_{22}K\Acl & P_{22}KB  
\end{bmatrix}
\begin{bmatrix}
    x\\
    \dz
\end{bmatrix}
\label{eq:Vdot}
\end{align}
where 
%we selected $Q = - \frac{1}{2}\left(P \smallmat{\Acl & B \\ 0 & I} + \smallmat{\Acl & B \\ 0 & I}^\top P\right)$ and 
we used the shortcut notation $\dz$ in place of $\dz(Kx)$.
With function $V$ in \eqref{eq:Lyapunov_f}, we may state our following second main result.

\begin{theo} \label{th:GAS}
The origin is globally asymptotically stable (GAS) for system \eqref{eq:linearCL} only if $\Acl$ is Hurwitz. Moreover, the origin is GAS if there exists a function $V$ as in \eqref{eq:Lyapunov_f} satisfying 
item (iv) of Theorem~\ref{theo:definite}, such that the opposite of its directional derivative $-
\dot V(x)$ in \eqref{eq:Vdot} satisfies item (ii) of Theorem~\ref{theo:definite}.
\end{theo}

\begin{proof}
The necessity that $\Acl$ Hurwitz follows from the fact that $\sat(u)=u$ in a neighborhood of the origin, therefore local asymptotic stability holds if and only if the (locally) linear dynamics governed by $\Acl$ is exponentially stable.

Let us now prove the sufficiency. First of all, by item (iv) of Theorem~\ref{theo:definite}, $V$ is positive definite and radially unbounded. 
Moreover, due to item (ii) of Theorem~\ref{theo:definite} $\dot V(x)$ 
satisfies:
\begin{align}
    \dot V(x) \leq -\lambda_{\min}(R)|\mathrm{sat}(Kx)|^2,\mbox{ for almost all }x\in \real^n.
\label{eq:Vdot-almost}
\end{align}
Due to \cite[Prop. 1]{DellaRossaMCSS21} and the continuity of the right-hand side of \eqref{eq:linearCL}, bound \eqref{eq:Vdot-almost} holds everywhere in the sense of Clarke. We may then apply the (continuous-time special case of the) non-smooth invariance principle in \cite[Thm 1]{SeuretTAC19} to conclude GAS. In fact, 
non-positivity of $\dot V$ implies that the origin is Lyapunov stable and that all trajectories are bounded and converge towards the largest invariant set $\mathcal{S}$ contained in $\{x\in\mathbb{R}^n: \sat(Kx)=0\}= \{x\in\mathbb{R}^n: Kx=0\}$, which means that the solution behaves linearly in $\mathcal{S}$.
Due to the fact that $\Acl=A_0-BK$ is Hurwitz, then the pair $(K,A_0)$ is detectable and, as a consequence $\mathcal{S}$ must be equal to the origin, thus proving global convergence, which together with Lyapunov stability implies GAS.
%
% Let us now characterize the set $\mathcal{S}$. As a direct consequence of the established inequality, we have $\dot{V}(x)=0\Rightarrow dz^\top W\ \mathrm{sat}(Kx)=0$ (since the term provided by the sector condition is non-negative by definition). This means in particular that
% $$
% -x^\top Q_0 x-\mathrm{sat}(Kx)^\top(2W+U)\,\mathrm{sat}(Kx)=0 
% $$
% which implies, in turn, that $\mathrm{sat}(Kx)=0$ for $x\in\mathcal{S}$ due to assumption \eqref{eq:cond3}. On the other hand, $\mathrm{sat}(Kx)=0$ is equivalent to $Kx=0$, so that \eqref{eq:sysmodel} necessarily behaves linearly, {\it i.e.} $dz\equiv 0$, within the set $\mathcal{S}$. The pair $(K,A_0)$ is detectable by construction\footnote{In particular $A_0-BK$ is Hurwitz}, so that $\{0\}$ is the only invariant set contained in $\mathcal{S}$, thereby proving the global asymptotic stability of \eqref{eq:sysmodel}. This concludes the proof.
\end{proof}

\section{Single-input results}
\label{sec:SI}

We discuss in this section a number 
of well-known single input relevant cases (single and double integrator, oscillator, integral oscillator) where we manage giving necessary and sufficient conditions for global asymptotic stability.

\subsection{Lyapunov construction for single-input systems}

To construct a Lyapunov function, we leverage condition \eqref{eq:LiuAkasaka} discussed in the introduction and 
we exploit Theorem~\ref{th:GAS} with the selection
$P_{11} = P_0+h K^\top K$ and $P_{22} = -h$ in \eqref{eq:Lyapunov_f}, where $h\in \real$, $P_0$ satisfies \eqref{eq:LiuAkasaka} and $P_0+hK^\top K\succ 0$, leading to
\begin{align}
    V(x) = \frac{1}{2}\left( x^\top P_0x + h(x^\top K^\top  K x-\dz(Kx)^2)\right).
\label{eq:Vspecific}
\end{align}
In the single input case, $2\int_0^{Kx}\sat(s)ds = |Kx|^2 - \dz(Kx)^2$, thus
selection \eqref{eq:Vspecific} is a convenient rewriting of the Popov-induced Lyapunov function in \cite[p. 277]{Khalil3} (see also the discussion in \cite[\S III.B]{ValmorbidaTAC22}), but we emphasize that the results of \cite[p. 277]{Khalil3} insist on positive definiteness of $P_0$, which we do not require here, thus making our tools applicable with ANCBI systems that cannot be globally exponentially stabilized. The relevance of using the expression in 
\eqref{eq:Vspecific} is in the possibility of applying the results of Section~\ref{sec:nonquadratic_lemmas} and their generalized form in Section~\ref{sec:generalization}, which leads to convenient LMI-based stability conditions.

 The structure in \eqref{eq:Vspecific} provides several advantages, the first one concerning the positive definiteness and radial unboundedness property, clarified next.

\begin{lemma}
\label{lem:pos_def}
    Function $V$ in \eqref{eq:Vspecific} is positive definite and radially unbounded if $P_0 \succeq 0$ and $P_0+hK^\top K\succ 0$.
\end{lemma}

\begin{proof}
When $h\leq  0$, we have $2V(x) = x^\top (P_0+hK^\top K)x + |h||\dz(Kx)|^2 \geq x^\top (P_0+hK^\top K)x$, which is trivially positive definite and radially unbounded.
When $h> 0$,
the proof is an immediate consequence of item (iv) of Theorem~\ref{theo:definite} with the selections 
$Q = \tfrac{1}{2}\smallmat{ P_0+hK^\top K & 0 \\ 0 & -h }$ and $T_0 = \tfrac{h}{2}>0$, which provides $Q - \Sigma_0 = \frac{1}{2}\smallmat{P_0 & 0\\ 0 & 0} + \frac{h}{2}\smallmat{K^\top \\ -1}\smallmat{K & -1}\succeq 0$. 
\end{proof}

A second advantage of the function in \eqref{eq:Lyapunov_f} concerns its directional derivative $\dot V$ along \eqref{eq:linearCL}. 
We prove in the proposition below, that,
when imposing the following {\em matching condition}
\begin{equation}\label{eq:match0}
\boxed{P_{0}B- hA_0 ^\top K^\top=\mu K^\top, \quad \mu \geq 0, }
\end{equation}
and by introduction the following notation and constraint
\begin{align}
\label{eq:omega_def}
    \boxed{\omega := 2 hKB,  \qquad \qquad  2\mu+\omega>0,}
\end{align}
the directional derivative, corresponding to \eqref{eq:Vdot} evaluated with the selections in \eqref{eq:Vspecific}, reads as
\begin{align}
\label{eq:Vdot_matching}
    \dot V(x) &= -\tfrac{1}{2}\smallmat{x \\ \dz(Kx)} Q \smallmat{x \\ \dz(Kx)},\mbox{ with}\\ 
    \nonumber
     &\; Q =  \begin{bmatrix}
    S_0 + (2\mu+\omega)K^\top K & -(\mu+\omega)K^\top\\ 
     -(\mu+\omega)K   & \omega
    \end{bmatrix},
\end{align}
where $S_0\succeq 0$ is defined in \eqref{eq:LiuAkasaka}.

The next lemma, in addition to proving \eqref{eq:Vdot_matching}, also provides conditions on $\mu$ and $\omega$ ensuring suitable decrease properties of $V$
by exploiting inequality \eqref{eq:LiuAkasaka}.

\begin{lemma}
\label{lem:VdotSpecific}
Under the matching condition \eqref{eq:match0}, \eqref{eq:omega_def}, the directional derivative of \eqref{eq:Vspecific} along \eqref{eq:linearCL} simplifies to \eqref{eq:Vdot_matching}. Moreover, if %2\mu+\omega>0$ and 
\eqref{eq:LiuAkasaka} holds, then $-
\dot V(x)$ in \eqref{eq:Vdot_matching} satisfies item (ii) of Theorem~\ref{theo:definite}.
\end{lemma}

\begin{proof}
Let us first prove \eqref{eq:Vdot_matching} by expanding \eqref{eq:Vdot} with 
$P_{11} = P_0+h K^\top K$ and $P_{22} = -h$, as in \eqref{eq:Vspecific}, under the matching condition \eqref{eq:match0}, which yields $-\dot V(x) = 
\frac{1}{2}\smallmat{x \\ \dz}^\top \smallmat{Q_{11} & Q_{12} \\ Q_{12}^\top & Q_{22}}  \smallmat{x \\ \dz}$, with
%~\footnote{We denote $\He(X)= X+X^\top$.}
\begin{align*}
    Q_{11} &= \He (-P_{11}\Acl) = \He ( -(P_0+hK^\top K)(A_0-BK))\\
    &= \He( -P_0A_0+P_0BK -h A_0^\top \KT K + \tfrac{\omega}{2}\KT K)\\ 
    Q_{12}&= -(P_0 + h \KT K) B + h (A_0-BK)^\top \KT \\ 
    &= -P_0B +hA_0^\top \KT -\omega \KT \\ 
    Q_{22} &= \omega.
\end{align*}
Substituting the matching condition \eqref{eq:match0} in the identities above leads to \eqref{eq:Vdot_matching}.

To prove the second part of the proposition, consider again the expression of $Q$ in \eqref{eq:Vdot_matching}. Exploiting the fact that $S_0\geq 0$, we may write
\begin{align*}
Q &\geq     \begin{bmatrix}
  (2\mu+\omega)K^\top K & -(\mu+\omega)K^\top\\ 
     -(\mu+\omega)K   & \omega
    \end{bmatrix} \\ 
    &= (2\mu+\omega) \bigmat{\KT K & \KT \\ K & 1} + \mu \bigmat{0 & \KT \\ K  & -2},
\end{align*}
which,
under the stated conditions $2\mu+\omega>0$ and $\mu\geq 0$ from
\eqref{eq:match0}, \eqref{eq:omega_def} implies that
item (ii) of Theorem~\ref{theo:definite} holds with $R = 2\mu+\omega>0$ and $T_0=\mu \geq 0$, as to be proven.
\end{proof}

\begin{rem}\label{rem:mu=0}
Let us stress that, whenever $P_0$ is not full rank, the matching condition \eqref{eq:match0} may only be satisfied with $\mu=0$. Indeed, 
we can notice that the inclusion of kernels in \eqref{eq:LiuAkasaka} implies that one can always find a matrix $M_0\in\mathbb{R}^{n\times n}$ satisfying the identity $A_0=M_0 P_0$. Using the factorization $M_0 P_0$ of $A_0$ in the matching condition \eqref{eq:match0} entails
$
P_0(B-hM_0^\top K^\top)=\mu K^\top
$, which in turn implies that $\mu K^\top\in\mathrm{Im}(P_0)$. On the other hand, when $P_0$ is not full rank, the latter condition and $P_0+h K^\top K\succ 0$ are incompatible, unless $\mu=0$.
    % \red{Comment on $\mu=0$ when $P_0$ is not full rank... (to be recalled in the examples). Either we can explain the ``ker'' condition in \eqref{eq:LiuAkasaka}, or we should remove it from \eqref{eq:LiuAkasaka}}
\end{rem}

In the remainder of this section, we illustrate the usefulness of the structure in \eqref{eq:Vspecific}
on a number of popular single-input examples that cannot be globally exponentially stabilized. We use Theorem~\ref{th:GAS} and Lemmas~\ref{lem:pos_def} and~\ref{lem:VdotSpecific} for each one of the considered cases.

\subsection{Saturated integrator}

As a first simple example consider the {\em saturated single integrator} 
$\dot x = \sat(-kx)$, with $x\in \real$, namely \eqref{eq:linearCL} with
$$
A_0=0,\quad B=1,\quad K=k, \quad \Acl = -k,
$$
briefly discussed in the introduction. For this example, the necessary condition $\Acl$ Hurwitz holds if and only if $k>0$. In fact, we may prove GAS for any $k>0$, which yields a necessary and sufficient condition. Rather than using the simple quadratic $V(x)=x^2$, let us use the arguably more elegant function
\begin{align}
\nonumber
 V(x) &=    \int_0^{kx}\!\!\! \sat(s)ds = \int_0^{kx} \!\!\!s - \dz(s)ds\\
 &= \tfrac{1}{2}((kx)^2 - \dz^2(kx)), 
 \label{eq:single_int_V}
\end{align}
which corresponds to \eqref{eq:Vspecific} with $P_0=0$ and $h=1$.
Then Lemma~\ref{lem:pos_def} applies because $P_0+h\KT K = k^2>0$. Moreover,
\eqref{eq:LiuAkasaka} holds with $S_0 = 0$ and
the matching condition \eqref{eq:match0} trivially holds with $\mu=0$ and we may evaluate
from \eqref{eq:omega_def}
 $\omega= 2k>0$. Since $2\mu + \omega = \omega>0$, then Lemma~\ref{lem:VdotSpecific} implies that $V$ in \eqref{eq:single_int_V} satisfies also the second assumption of
 Theorem~\ref{th:GAS} and GAS follows.

%  Then item (iv) of Theorem~\ref{theo:definite} establishes positive definiteness and radial unboundedness of $V$ with $T_0=1$.
 
% The directional derivative of $V$ along $\dot x = -\sat(kx)$ reads
% \begin{align}
% \nonumber
%  \dot V(x) &= \sat(kx)k\dot x = -k|\sat(kx)|^2 = -k(kx-\dz(kx))^2\\
%       &= - k \smallmat{x \\ \dz(kx)}^\top\smallmat{k \\ -1}\smallmat{k -1} \smallmat{x \\ \dz(kx)}
% \end{align}
% which corresponds to \eqref{eq:Vdot} with $Q = \smallmat{k^2 & -k\\ -k 
%  & 1}$. Then item (iii) of Theorem~\ref{theo:definite} establishes positive definiteness of $-\dot V$ with $R=1$. GAS finally follows from Theorem~\ref{th:GAS}.

 \subsection{Saturated double integrator}
\label{sec:double_int}

Consider now a saturated double integrator, namely \eqref{eq:linearCL} with
\begin{align}
\label{eq:double_int}
A_0=\begin{bmatrix}
0\! & 1\\
0\! & 0
\end{bmatrix},\; B=\begin{bmatrix}0\\ 1\end{bmatrix},\; K=[k_1\ k_2],\; \Acl = \begin{bmatrix}
0 & 1\\
-k_1 \!\!&\!\! -k_2
\end{bmatrix}
\end{align}
Using Routh criterion, it is immediate to check that $\Acl$ is Hurwitz if and only if $k_1>0$ and $k_2>0$. In fact, the following is a well known fact.

\begin{prop}
\label{prop:doubleInt}
The saturated double integrator \eqref{eq:linearCL}, \eqref{eq:double_int} is GAS if and only if $k_1>0$ and $k_2>0$.
\end{prop}

The well-known proof of Proposition~\ref{prop:doubleInt} (see, e.g., \cite{ShiSaberiACC02}) is based on the Popov-like Lyapunov function
\begin{align}
    V(x) =  \frac{x_2^2}{2} + \int_0^{k_1x_1+k_2x_2}\hspace{-.6cm}\sat(s)ds,
\end{align}
for which, using similar tricks to those in \eqref{eq:single_int_V}, we can obtain the expression
$V(x)=\frac{1}{2}\left( x_2^2 +\frac{1}{k_1}(x^\top K^\top  K x-\dz(Kx)^2)\right)$,
which corresponds to \eqref{eq:Vspecific} with the parameters
% with the selection $h=\displaystyle\tfrac{1}{k_1}$. It is straightforward to check that $V$ in \eqref{eq:Vdouble_int_ours} corresponds to \eqref{eq:Lyapunov_f}
% with 
\begin{align}
% P= \begin{bmatrix}
% P_{11}&0\\
% 0&-\frac{1}{2k_1} 
% \end{bmatrix},
% \mbox{ and } 
P_0=\begin{bmatrix}
0 & 0\\
0 & 1
\end{bmatrix}, \quad 
h=\displaystyle\frac{1}{k_1}.
\label{eq:Vdouble_int_ours}
\end{align}
%
% , or equivalently in matrix form
% $$
% V(x)=[x^\top\ dz]\underbrace{\begin{bmatrix}
% P_{11}&0\\
% 0&-\frac{1}{2k_1} 
% \end{bmatrix}}_{=:\myQ }\begin{bmatrix}
% x\\
% dz
% \end{bmatrix}
% $$
% with 
% $$
% P_{11}=\frac1{2k_1}\begin{bmatrix}
% k_1^2 & k_1k_2\\
% k_2k_2&k_2^2+k_1
% \end{bmatrix}
% $$
We see that $P_{11}=P_0 + hK^\top K\succ 0$, so that Lemma~\ref{lem:pos_def} establishes
% Application of item (iv) of Theorem~\ref{theo:definite} with $T_0=h/2=(2k_1)^{-1}$ shows that $V$ is 
positive definiteness and radial unboundedness of~$V$. 

Also for this example we obtain $S_0 = \smallmat{0 & 0 \\ 0 & 0}$ so that  
\eqref{eq:LiuAkasaka} holds. Moreover 
the matching condition~\eqref{eq:match0} holds with $\mu = 0$ again, consistently with  Remark~\ref{rem:mu=0}. Computing $\omega$ from \eqref{eq:omega_def}, we obtain $\omega = 2\frac{k_2}{k_1}>0$, which is positive under the necessary conditions for $\Acl$ Hurwitz. Then Lemma~\ref{lem:VdotSpecific} establishes once again that the second assumption of Theorem~\ref{th:GAS} holds too and GAS follows.

% It is interesting (and important for the generalizations developed in Section~\ref{sec:generalization}) to check that
% \begin{equation}\label{eq:Lyap_eq}
% \boxed{P_{11}\Acl +\Acl ^\top P_{11}=-2r K^\top K}
% \end{equation}
% with
% $r=\displaystyle\frac{k_2}{k_1}$.
% Evaluating $\dot{V}(x)$, we obtain the expression in \eqref{eq:Vdot} with
% \begin{align}
% Q=\begin{bmatrix}
% P_{11}\Acl +\Acl ^\top P_{11}&  P_{11}B- h \Acl ^\top K^\top\\
% B^\top P_{11} - h K \Acl & -2hKB
% \end{bmatrix}.
% \label{eq:Qdouble_int}
% \end{align}
% % where we used the shortcut $\dz$ instead of $\dz(Kx)$.
% Checking that the {\it matching} identity 
% \begin{equation}\label{eq:match}
% \boxed{P_{11}B- h\Acl ^\top K^\top=\lambda K^\top}
% \end{equation}
% holds with $\lambda=2r=2k_2/k_1>0$, and observing that $hKB=r$, we get
% $$
% Q=\lambda \begin{bmatrix}
% K^\top K& -K^\top\\
% -K & 1
% \end{bmatrix},
% $$
% thus showing that $-\dot{V}(x)$ satisfies item (ii) of Theorem~\ref{theo:definite} with $R = \lambda I$, $T_0=0$ and $T_1=0$. Then GAS follows from Theorem~\ref{th:GAS} with Proposition~\ref{prop:pos_def}.

\subsection{Saturated oscillator}
\label{sec:osc}

Consider the saturated oscillator having frequency $\omega>0$, corresponding to
\eqref{eq:linearCL} with
% $$
% \dot{x}=\underbrace{(A_0 -BK)}_{=:\Acl }x+B dz
% $$
% with
\begin{align}
A_0 =\begin{bmatrix}
0 & \omega \\
-\omega & 0
\end{bmatrix},\; B=\begin{bmatrix}0\\ \omega \end{bmatrix},\; 
\Acl = \begin{bmatrix}
0 & \omega \\
-k_1\omega -1 & -k_2\omega
\end{bmatrix}
\label{eq:oscillatoromega}
\end{align}
with the usual notation $K=[k_1\ k_2]$ and where we note that the input gain $\omega$ can be adjusted at will, due to the arbitrariness of the saturation limits.
As a first step, we notice that we may perform a time-scale change $\tau = \omega t$ (so that $x' = \frac{dx}{d\tau} = \frac{1}{\omega} \frac{dx}{dt} = \frac{1}{\omega} \dot x$) to transform the system into the normalized dynamics corresponding to
\eqref{eq:linearCL} with
\begin{align}
A_0 =\begin{bmatrix}
0 & 1 \\
-1 & 0
\end{bmatrix},\; B=\begin{bmatrix}0\\ 1 \end{bmatrix},\; 
\Acl = \begin{bmatrix}
0 & 1 \\
-k_1 -1 & -k_2
\end{bmatrix}.
\label{eq:oscillator}
\end{align}
We may easily check that $\Acl$ is Hurwitz if and only if $k_1>-1,\ k_2>0$.
Also for this case we may prove the following proposition, paralleling Proposition~\ref{prop:doubleInt}.

\begin{prop}
\label{prop:oscill}
The saturated oscillator \eqref{eq:linearCL}, \eqref{eq:oscillatoromega} is GAS if and only if $k_1>-1$ and $k_2>0$.
\end{prop}

\begin{proof}
To prove Proposition~\ref{prop:oscill} we construct $V$ as in \eqref{eq:Vspecific} with
\begin{align}
\label{eq:Voscill}
P_0 = I, \quad h = \frac{k_1}{|K|^2},
% V(x)=\frac{x^\top x}{2} +\frac{h}{2}(x^\top K^\top  K x-\dz(Kx)^2)
\end{align}
which is positive definite and radially unbounded due to Lemma~\ref{lem:pos_def}, because $P_0\succ 0$
and also $P_0 + h K^\top K = I + k_1\frac{K^\top K}{|K|^2}$ is positive definite for any selection of $K$
satisfying the necessary condition  $k_1>-1,\ k_2>0$. Indeed, matrix $\frac{K^\top K}{|K|^2}$ is a projection with one eigenvalue at zero and one eigenvalue at 1 so that $k_1>-1$ implies  $k_1\frac{K^\top K}{|K|^2} \succ -I$.
We also emphasize that the sign of $h$ corresponds to the sign of $k_1$.

% with the selection $h=\displaystyle\frac{k_1}{k_1^2+k_2^2}$. Once again, $V$ in \eqref{eq:Voscill} corresponds to \eqref{eq:Lyapunov_f}
% with 
% \begin{align*}
% P= \begin{bmatrix}
% P_{11}&0\\
% 0& -\frac{h}{2} 
% \end{bmatrix},
% \mbox{ and } 
% P_{11}=\frac1{2}\begin{bmatrix}
% 1+hk_1^2 & hk_1k_2\\
% hk_2k_2&1+hk_2^2
% \end{bmatrix}
% \end{align*}
% Application of item (iv) of Theorem~\ref{theo:definite} with $T_0= ???? h/2=(2k_1)^{-1}$ shows that $V$ is positive definite and radially unbounded. 

Let us now focus on $\dot V$.
Also for this example we obtain $S_0 = \smallmat{0 & 0 \\ 0 & 0}$ so that  
\eqref{eq:LiuAkasaka} holds. 
After some calculations, we find that the matching condition
\eqref{eq:match0} holds with $\mu = \frac{k_2}{|K|^2}$. We may also easily compute from \eqref{eq:omega_def} $\omega = 2\frac{k_1k_2}{|K|^2}$, which shows that the condition
$2\mu + \omega = \frac{2(1+k_1)k_2}{|K|^2}>0$ holds true under the necessary conditions for $\Acl$ being Hurwitz, given in the statement. 
Since all of the conditions of Lemma~\ref{lem:VdotSpecific} hold, then all of the assumptions of Theorem~\ref{th:GAS} are satisfied and GAS follows.
\end{proof}

We emphasize that Proposition~\ref{prop:oscill} provides necessary and sufficient global asymptotic stability conditions for the saturated oscillator with linear feedback. To the best of our knowledge, this necessary and sufficient result is new and relies on a Lyapunov function of the form \eqref{eq:Vspecific}
with $P_0$ satisfying \eqref{eq:LiuAkasaka}.

% \red{

% The candidate Lyapunov function is
% $$
% V(x)=\frac{x^\top x}{2} +\frac{h}{2}(x^\top K^\top  K x-dz^2)
% $$
% with $h=\displaystyle\frac{k_1}{k_1^2+k_2^2}$, or equivalently in matrix form
% $$
% V(x)=[x^\top\ dz]\underbrace{\begin{bmatrix}
% P_{11}&0\\
% 0&-\frac{h}{2} 
% \end{bmatrix}}_{=:\myQ }\begin{bmatrix}
% x\\
% dz
% \end{bmatrix}
% $$
% with 
% $$
% P_{11}=\frac1{2}\begin{bmatrix}
% 1+hk_1^2 & hk_1k_2\\
% hk_2k_2&1+hk_2^2
% \end{bmatrix}
% $$
% By direct inspection, we find again a condition like \eqref{eq:Lyap_eq}:
% $$
% P_{11}\Acl +\Acl ^\top P_{11}=-rK^\top K,
% $$
% now with $r=\displaystyle\frac{(1+k_1)k_2}{k_1^2+k_2^2}=\textcolor{teal}{\frac{k_2}{k_1^2+k_2^2}+hk_2}$. Moreover, in the evaluation of $\dot{V}(x)$, the matching condition \eqref{eq:match} holds again in this case:
% $$
% P_{11}B-\frac{h}{2}\Acl ^\top K^\top=\lambda K^\top
% $$
% with
% $$
% \lambda=\frac{(1+2k_1)k_2}{2(k_1^2+k_2^2)}=\textcolor{teal}{\frac{k_2}{2(k_1^2+k_2^2)}+hk_2}
% $$
% This implies that 
% $$ 
% \dot{V}(x)=[x^\top\ dz]\begin{bmatrix}
% -r K^\top K&\lambda K^\top\\
% \star&-hk_2 
% \end{bmatrix}\begin{bmatrix}
% x\\
% dz
% \end{bmatrix}
% $$
% \textcolor{red}{Caso $k_1>0$} By invoking an upper sector condition \textcolor{teal}{with $T_1=\frac{k_2}{2(k_1^2+k_2^2)}$}, proves that $\dot{V}(x)$ is negative semidefinite.\bigskip\\
% \textcolor{red}{Caso $-1<k_1\leq0$} By invoking a lower sector condition \textcolor{teal}{with $T_0=\frac{k_2}{2(k_1^2+k_2^2)}$}, proves that $\dot{V}(x)$ is negative semidefinite.
% }

% \newpage
% dasdaa

% \newpage

\subsection{Saturated integral oscillator}

Let us now consider the saturated integral oscillator comprising the cascaded of an integrator and an oscillator, corresponding to \eqref{eq:linearCL} with
\begin{align}
\label{eq:integroscillatoromega}
[A_0|B] &=\left[ \begin{array}{ccc|c}
0 & \omega & 0 & 0\\
\!\!\!-\omega & 0 & \zeta  & 0\\
0&0&0 & \frac{\zeta}{\omega^2}\!\!
\end{array}
\right],\\ %B=\begin{bmatrix}0\\0\\ \frac{\zeta}{\omega^2}\end{bmatrix},\;
\Acl &=\begin{bmatrix}
0 & \omega & 0 \\
-\omega & 0 & \zeta \\
-\frac{\zeta}{\omega^2}k_1& \!\!\!-\frac{\zeta}{\omega^2}k_2& \!\!\!-\frac{\zeta}{\omega^2}k_3
\end{bmatrix},
\nonumber
\end{align}
where $\omega>0$ and $\zeta \neq 0$ are arbitrary quantities,
and with the usual notation $K=[k_1\; k_2\; k_3]$. 
Also in this case the input gain $\frac{\zeta}{\omega^2}$ does not affect the analysis modulo a change in the saturation limits.

Performing a
 time-scale change $\tau = \omega t$, as in 
\eqref{eq:oscillatoromega},
and with the state transformation $\xi_1=x_1$, $\xi_2=x_2$, $\xi_3 = \frac{\zeta}{\omega}x_3$,
we may again transform system \eqref{eq:integroscillatoromega} into a normalized form corresponding to \eqref{eq:linearCL} with
\begin{align}
A_0 =\begin{bmatrix}
0 & 1 & 0 \\
-1 & 0 & 1\\
0&0&0
\end{bmatrix},\; B=\begin{bmatrix}0\\0\\ 1\end{bmatrix},\; 
\Acl =\begin{bmatrix}
0 & 1 & 0 \\
-1 & 0 & 1\\
-k_1 & -k_2 & -k_3
\end{bmatrix}.
\label{eq:integroscillator}
\end{align}
Using Routh criterion, the necessary condition that $\Acl$ be Hurwitz, corresponds to requiring the following three conditions on the gain $K=[k_1\; k_2\; k_3]$:
\begin{align}
\label{eq:satosci_Kconds}
\mbox{(i) } k_3>0, \quad
\mbox{(ii) } k_2k_3>k_1, \quad
\mbox{(iii) } k_1+k_3>0.
\end{align}
For our construction, we also require the following condition
\begin{align}
\label{eq:satosci_ExtraK1cond}
\mbox{(iv*) } k_1\leq 0,
\end{align}
which can be restrictive whenever it holds that $k_2k_3>0$. In fact, the necessary conditions \eqref{eq:satosci_Kconds} imply 
$-k_3 < k_1<k_2 k_3$, while the additional condition \eqref{eq:satosci_ExtraK1cond} rules out possible positive selections of $k_1$ in this range.
Also for this system our Lyapunov construction provides powerful stability certificates, as proven next.

\begin{prop}
\label{prop:integroscill}
The saturated integral oscillator \eqref{eq:linearCL}, \eqref{eq:integroscillatoromega} is GAS if $K$ satisfies \eqref{eq:satosci_Kconds}, \eqref{eq:satosci_ExtraK1cond}.
\end{prop}

\begin{proof}
To prove Proposition~\ref{prop:integroscill},
due to the equivalence between \eqref{eq:integroscillatoromega} and \eqref{eq:integroscillator}, we focus on the
normalized dynamics \eqref{eq:integroscillator} and
we construct $V$ as in \eqref{eq:Vspecific} with
\begin{align}
\label{eq:Vintegroscill}
P_0 = \begin{bmatrix}
1 & 0 & -1\\
0 & 1 & 0\\
-1 & 0 & p_{33}
\end{bmatrix},\; p_{33}=\displaystyle\frac{k_2^2-k_1k_3}{k_1^2+k_2^2},\;
h=\displaystyle\frac{k_2}{k_1^2+k_2^2},
\end{align}
where we emphasize that $S_0$ in \eqref{eq:LiuAkasaka} is zero for any selection of $p_{33}$. However, the selected value in \eqref{eq:Vintegroscill} allows proving the additional properties of $V$ via Lemmas~\ref{lem:pos_def} and~\ref{lem:VdotSpecific}, as discussed next.

Let us start by studying the positive definiteness of $V$ in \eqref{eq:Vspecific} with the selections \eqref{eq:Vintegroscill}
exploiting Lemma~\ref{lem:pos_def}. First note that condition (iii) in \eqref{eq:satosci_Kconds} and \eqref{eq:satosci_ExtraK1cond} imply
$$
|k_1|=-k_1<k_3 = |k_3| \quad \Rightarrow \quad -k_1k_3 = |k_1||k_3| \geq |k_1|^2,
$$
which implies $p_{33}\geq 1$ and, as a consequence, $P_0 \succeq 0$ (the first assumption of Lemma~\ref{lem:pos_def}). To prove the other assumption (namely $P_0 +h K^\top K\succ 0$),
we first observe that due to (ii) in \eqref{eq:satosci_Kconds}, $k_1$ and $k_2$ 
cannot be both zero, so it is enough to prove
positive definiteness of $(k_1^2+k_2^2)(P_0 + hK^\top K)$. Let us also note that, due to \eqref{eq:satosci_Kconds}, we have
\begin{align}
    k_2k_3 > k_1>-k_3 \quad \Rightarrow \quad 1+k_2>0.
\label{eq:k2-1positive}
\end{align}
Then, we apply Sylvester criterion to matrix
$(k_1^2+k_2^2)(P_0 + hK^\top K)$,
whose three principal minors, after some cumbersome calculations, can be expressed as
\begin{align*}
m_1 &= k_2^2 + k_1^2(1+k_2), \\
m_2 &= (k_1^2 +k_2^2)^2(1+k_2),\\
m_3 &=  (k_1^2 +k_2^2)^2(k_1+k_3)(k_2k_3-k_1) ,
\end{align*}
which are easily shown to be all positive under the three conditions in  \eqref{eq:satosci_Kconds} and the observation in \eqref{eq:k2-1positive}, together with the fact that $k_1$ and $k_2$ cannot be both zero from (ii) in \eqref{eq:satosci_Kconds}.
Consequently, all the assumptions of Lemma~\ref{lem:pos_def} hold, so that $V$ is positive definite and radially unbounded.

Let us now focus on $\dot V$ and recall that 
$S_0 = \smallmat{0 & 0 & 0 \\ 0 & 0& 0 \\ 0 & 0& 0}$ , which implies \eqref{eq:LiuAkasaka}.
After some calculations, we find that the matching condition
\eqref{eq:match0} holds
with the selection
$\mu = -\frac{k_1}{k_1^2+k_2^2} \geq 0$, which is non-negative because of condition~\eqref{eq:satosci_ExtraK1cond}.
We may also compute, from \eqref{eq:omega_def}, $\omega = \frac{2k_2k_3}{k_1^2+k_2^2}$
and we note that the condition in Proposition~\ref{lem:VdotSpecific},
$2\mu + \omega = 2\frac{k_2k_3-k_1}{k_1^2+k_2^2}>0$ holds true because of item (ii) in \eqref{eq:satosci_Kconds}. As a consequence, Lemma~\ref{lem:VdotSpecific} implies also the second assumption of Theorem~\ref{th:GAS}, and GAS is proven.
%
% s \eqref{eq:Lyap_eq} and \eqref{eq:match} hold also in this case, with the scalars
% \begin{align}
% \label{eq:rlambda_integroscillator}
% r=\frac{k_2k_3-k_1}{k_1^2+k_2^2}>0, \;
% \lambda=\frac{2k_2k_3-k_1}{k_1^2+k_2^2} = r+ \frac{k_2k_3}{k_1^2+k_2^2},
% % \label{eq:lambda_integroscillator}
% \end{align}
% where positivity of $r$ follows from item (ii) of \eqref{eq:satosci_Kconds}.
% Exploiting the matching conditions \eqref{eq:Lyap_eq} and \eqref{eq:match}, it is possible to express the directional derivative 
% $\dot V$ 
% of $V$ as in \eqref{eq:Vdot}, with $Q$ as in \eqref{eq:Qdouble_int}, which, due to 
% the matching conditions \eqref{eq:Lyap_eq} and \eqref{eq:match} with $r$ and $\lambda$ as in \eqref{eq:rlambda_integroscillator}, and due to
% $hKB = \frac{k_2 k_3}{k_1^2 + k_2^2}$, yields
% \begin{align*}
%     Q = \frac{1}{k_1^2 + k_2^2}\begin{bmatrix}
%         2(k_2k_3-k_1)K^\top K & -(2k_2k_3-k_1)K^\top \\
%         (-k_1+2k_2k_3)K & 2k_2k_3
%     \end{bmatrix}
% \end{align*}
\end{proof}

% Imposing the matrix inequalities \eqref{eq:cond1}-\eqref{eq:cond5}, Theorem~\ref{theo:stability} provides the candidate Lyapunov function
% $$
% V(x)=x^\top P_0 x +(x^\top K^\top H K x-dz^\top H dz)=[x^\top\ dz]\underbrace{\begin{bmatrix}
% P_{11}&0\\
% 0&-H 
% \end{bmatrix}}_{=:\myQ }\begin{bmatrix}
% x\\
% dz
% \end{bmatrix}
% $$
% with 
% $$
% P_0:=\begin{bmatrix}
% 1 & 0 & -1\\
% 0 & 1 & 0\\
% -1 & 0 & p_{33}
% \end{bmatrix},\quad P_{11}=P_0+ K^\top H K
% $$
% where $p_{33}=\displaystyle\frac{k_2^2-k_1k_3}{k_1^2+k_2^2}$, $H=\displaystyle\frac{k_2}{k_1^2+k_2^2}$ and $W=\displaystyle\frac{-k_1}{k_1^2+k_2^2}$.\smallskip\\ It is worth stressing that, under condition iv*), one has $w\geq0$ and $p_{33}\geq 1$, which implies $P_0\succeq 0$. It is also interesting to note that $H\preceq 0$ for $k_2\leq0$: in this case, hinging on i)-iii), it can be easily checked that the necessary condition $P_{11}\succeq0$ still holds. Let us finally stress that, for any $p_{33}\in\mathbb{R}$, one has
%  $$
%  P_0 A_0 +A_0 ^\top P_0=0
%  $$
%  by construction. Overall, evaluating the derivative of ${V}(x)$ yields $\dot{V}(x)=\frac{1}{k_1^2+k_2^2}{\tilde{V}}(x)$ with
% $$
% \tilde{V}(x)=[x^\top\ dz]\begin{bmatrix}
% (2k_1-2k_2k_3)K^\top K&(-k_1+2k_2k_3)K^\top\\
% *&-2k_2k_3
% \end{bmatrix}\begin{bmatrix}
% x\\
% dz
% \end{bmatrix}
% $$
% Using the lower sector condition with $T_0=-k_1>0$ allows to conclude that $\tilde{V}(x)\leq 0 \Rightarrow \dot{V}(x)\leq 0$.

Proposition~\ref{prop:integroscill} is a new result, to the best of our knowledge. 
While it only provides sufficient conditions for GAS, necessity would hold if we could prove that assumption \eqref{eq:satosci_ExtraK1cond} is necessary for GAS of the saturated feedback. Note that \eqref{eq:satosci_ExtraK1cond} is not necessary for LES of the local linear dynamics,
thereby highlighting a peculiar case where saturation would jeopardize global convergence,
as commented in the next remark.
% We conjecture in the next remark that 
% Proposition~\ref{prop:integroscill} actually
% characterizes
% necessary and sufficient conditions.

% for GAS of \eqref{eq:linearCL}, \eqref{eq:integroscillatoromega}, it still covers most of the possible selections of the gain $K$,
% because assumption \eqref{eq:satosci_ExtraK1cond} is not necessary for 
% Some additional considerations about assumption \eqref{eq:satosci_ExtraK1cond} are given in the next remark.

\begin{rem} %[{\bf On the alleged necessity of condition} iv*)] 
[{\bf Conjecture about necessary and sufficient conditions}]
We conjecture that Proposition~\ref{prop:integroscill} actually provides necessary and sufficient conditions for GAS of the saturated linear feedback \eqref{eq:linearCL}, \eqref{eq:integroscillatoromega}. 
Our conjecture is motivated by the fact that when
% \red{qui direi piu' esplicitamente che cerchiamo necessary and sufficient GAS conditions -- for example ``Conjecture about necessary and sufficient conditions'' -- We conjecture that Proposition~\ref{prop:integroscill} actually provides necessary and sufficient conditions for GAS of the saturated linear feedback \eqref{eq:linearCL}, \eqref{eq:integroscillatoromega}. Such a result is motivated by the ????? -- magari citare Fuller dicendo che lui lo ha dimostrato per il triplo integratore che c'e' LAS ma non GAS con qualsiasi feedback lineare...}
% When 
$k_1>0$,
thus violating our assumption in \eqref{eq:satosci_ExtraK1cond},
 the matching condition \eqref{eq:match0} still holds, but with a negative
value of $\mu = -\frac{k_1}{k_1^2+k_2^2} < 0$, therefore Lemma~\ref{lem:VdotSpecific} cannot be applied
and we suspect that there exists a 
locally attractive oscillatory behavior induced by 
input oscillations repeatedly hitting the saturation limits.
 In fact,
all attempts to find a working sector condition failed for this case
%. Moreover,
and all of the numerical tests that we performed 
show the emergence of nonconverging responses for large enough initial conditions.
% (numerical) counterexamples to the stabilization of the saturated system could be found. 
This suggests that condition iv*) might in fact be necessary for bounded stabilization. It would be interesting to investigate deeper whether the lack of GAS 
can be proven by exploiting the fact that, when $k_1>0$, the matrix $P_0$ becomes sign-indefinite.
\end{rem}

% \begin{rem}
% The proof given above can be immediately extended to the case where the oscillator frequency is arbitrary, namely when 
% the open-loop matrix $A_0 $ in ??? is replaced by
% $A_{\omega_\circ} =\begin{bmatrix} 0 & \omega_\circ & 0 \\ -\omega_\circ & 0 & 1 \\ 0 & 0 & 0\end{bmatrix}$, because a time-scale change $\tau = \omega_\circ t$ (so that $x' = \frac{dx}{d\tau} = \frac{1}{\omega_\circ} \frac{dx}{dt} = \frac{1}{\omega_\circ} \dot x$) and the rescaled state variables $\xi_1 = \omega_\circ^2 x_1$, $\xi_2 = \omega_\circ^2 x_2$, $\xi_3 = \omega_\circ x_3$ transform the dynamics into $\xi' = A_0  \xi + Bu$.
% \end{rem}

\newcommand\Hpos{H_\mathrm{p}}
\newcommand\Hneg{H_\mathrm{n}}

\section{General architecture of Lyapunov functions}
\label{sec:generalization}

We generalize here the single-input construction given in Section~\ref{sec:SI}
and rely on the following intuitive multi-input generalization of the function in 
\eqref{eq:Vspecific},
\begin{align}
    V(x) = \frac{1}{2}\left( x^\top P_0x + x^\top K^\top H K x-\dz(Kx)^\top H\dz(Kx)\right),
\label{eq:VspecificMI}
\end{align}
where matrix $H$ is diagonal  and
matrix $P_0$ satisfies \eqref{eq:LiuAkasaka} and the following %\red{(relaxed version of \eqref{eq:LiuAkasaka} ?? } 
generalization of the conditions in  Lemma~\ref{lem:pos_def}.
\begin{align}
\label{eq:condPosMI1}
&P_0\succeq 0, \quad -S_0:=P_0A_0+A_0^\top P_0\preceq 0\\
\label{eq:condPosMI2}
& P_0+K^\top H K\succ 0,
\end{align}
where we note that \eqref{eq:condPosMI1} resembles the condition~\ref{eq:LiuAkasaka} commented in the introduction.
We prove below that Lemmas~\ref{lem:pos_def} and~\ref{lem:VdotSpecific} can be extended to the multi-input case by leveraging function $V$ in \eqref{eq:VspecificMI}.

To generalize Lemma~\ref{lem:pos_def}, %to the multi-input expression \eqref{eq:VspecificMI} 
we perform a decomposition of the diagonal matrix $H$
in a positive and negative components, satisfying the following:
\begin{align}
\label{eq:Hdecomposition}
H = \Hpos - \Hneg,\; \text{ with }\;
\begin{cases}
    \Hpos \succ 0, \; \Hneg \succeq 0, \\
    P_0-K^\top \Hneg K \succeq 0.
\end{cases}
\end{align}

\begin{rem}
While the decomposition 
$H = \Hpos - \Hneg$ always exists, in order to satisfy the second inequality in \eqref{eq:Hdecomposition}, 
one may attempt either to select $\Hpos$ as small as possible, hoping that inequality \eqref{eq:condPosMI2} can be exploited, or otherwise to select $\Hneg$
as small as possible, hoping that the left inequality in \eqref{eq:condPosMI1} can be exploited. Both strategies are effective, whenever feasible.
In general, the best strategy is to define two decision variables $\Hpos$, $\Hneg$ in a convex optimization formulation, as explained in Corollary~\ref{cor:LMI} at the end of this section.
\end{rem}

% \color{red}
% \begin{rem}
% Conditions \eqref{eq:condPosMI2}-\eqref{eq:H2} comprise the case of a diagonal $H\succ0$, as well as the case $H\preceq 0$. Indeed, in the first case, one may trivially select $H_1=H$ and $H_2=0$. Conversely, when $H\preceq 0$, due to continuity of \eqref{eq:condPosMI2}, a sufficiently small $\epsilon>0$ can always be found such that 
% $P_0+K^THK\succeq \epsilon K^TK$. Then the desired decomposition $H=H_1-H_2$ holds true with $H_1=\epsilon I$ and $H_2=-H+\epsilon I$.
% \end{rem}
% \color{black}

%
% can be extended to the multi-input case by leveraging function $V$ in \eqref{eq:VspecificMI}.
%
% Regarding Lemma~\ref{lem:pos_def}, 
 % it is enough to 
% We also impose the following constraints to ensure positive definiteness and radial unboundedness:
We can then state the following generalization of Lemma~\ref{lem:pos_def}.
% We state and proof the corresponding results next.

\begin{lemma}
    \label{lem:pos_def_general}
     If matrices $P_0$ and $H$ satisfy \eqref{eq:condPosMI2} and  \eqref{eq:Hdecomposition},
    % \eqref{eq:condPosMI1}, 
    then function $V$ in \eqref{eq:VspecificMI} satisfies 
    item (iv) of Theorem~\ref{theo:definite}. In particular, it
    is positive definite and radially unbounded. 
\end{lemma}

\begin{proof}
    Following parallel steps to those in the second part of proof of Lemma~\ref{lem:pos_def} and using the decomposition in \eqref{eq:Hdecomposition}, one has (for compact notation, we use the shortcut $\dz$ 
    in place of $\dz(Kx)$)
\begin{align*}
%\nonumber
    &2V(x)= %\frac12
    \begin{bmatrix}
    x\\
    \dz
\end{bmatrix}^\top 
\begin{bmatrix}
  P_0+K^\top (\Hpos - \Hneg) K & 0 \\
 0 & -\Hpos + \Hneg  
\end{bmatrix}
\begin{bmatrix}
    x\\
    \dz
\end{bmatrix}\\
%     &\geq %\frac12 
%     \begin{bmatrix}
%     x\\
%     \dz
% \end{bmatrix}^\top \begin{bmatrix}
%   P_0+K^\top (\Hpos - \Hneg)K & 0 \\
%  0 & -\Hpos 
% \end{bmatrix}
% \begin{bmatrix}
%     x\\
%     \dz
% \end{bmatrix}\\
    & \quad = %\frac12 
    \begin{bmatrix}
    x\\
    \dz
\end{bmatrix}^\top \!\!\!\!\left(\begin{bmatrix}
  P_0-K^\top\Hneg K\!\!\! & 0 \\
 0 &  \!\!\! \Hneg
\end{bmatrix} \!+\! 
\begin{bmatrix}
  K^\top\Hpos K \!\!\! & 0 \\
 0 &  \!\!\! -\Hpos
\end{bmatrix}
\right)
\begin{bmatrix}
    x\\
    \dz
\end{bmatrix},
\end{align*}
which fits the assumptions of 
% Now, hinging on \eqref{eq:H2}, 
item (iv) of Theorem~\ref{theo:definite} with the selections 
$Q=\smallmat{P_0+K^\top (\Hpos - \Hneg) K & 0 \\
 0 & -\Hpos + \Hneg}$ and $T_0=H_1$. Indeed,
$Q_{11}\succ 0$ from \eqref{eq:condPosMI2},
and
\begin{align*}
    Q-\Sigma_0 = \begin{bmatrix}
  P_0-K^\top\Hneg K & 0 \\
 0 &  \Hneg
\end{bmatrix} + 
\begin{bmatrix}
  K^\top\Hpos K & K^\top\Hpos \\
 \Hpos K &  -\Hpos
\end{bmatrix}
\end{align*}
 is positive semi-definite due to \eqref{eq:Hdecomposition}.
% $Q_{11}=P_0+K^\top (H_1-H_2)K$ and $T_0=H_1$, thereby proving that $V(x)$ is positive definite and radially unbounded.
\end{proof}

To generalize Lemma~\ref{lem:VdotSpecific} to the multi-input expression \eqref{eq:VspecificMI},
% Regarding Lemma~\ref{lem:VdotSpecific}, 
% the following generalization of 
we generalize the matching condition \eqref{eq:match0} as
 \begin{equation}
 \label{eq:match_MI}
P_0B-A_0^\top K^\top H =K^\top M,\ M\succeq 0,\; M\mbox{ diagonal}
\end{equation}
alongside with the following generalization of
the property $2\mu +\omega>0$ imposed in
\eqref{eq:omega_def} and required in 
Lemma~\ref{lem:VdotSpecific},
 \begin{equation}\label{eq:condOMEGA}
 2M+\underbrace{HKB+B^\top K^\top H}_{=:\Omega}\succ 0.
\end{equation}
The generalization, stated below, now involves the matrix $\Omega:=HKB+B^\top K^\top H$, generalizing the scalar $\omega$ in \eqref{eq:omega_def}. In particular, paralleling 
Lemma~\ref{lem:VdotSpecific}, we will show that \eqref{eq:Vdot_matching}
generalizes to
\begin{align}
\label{eq:Vdot_matchingMI}
    \dot V(x) &= -\tfrac{1}{2}\smallmat{x \\ \dz(Kx)} Q \smallmat{x \\ \dz(Kx)},\mbox{ with}\\ 
    \nonumber
     &\;
    Q=\begin{bmatrix}
        S_0 + K^\top (2M+\Omega)K &- K^\top (M+\Omega)\\
-(M+\Omega)K & \Omega
    \end{bmatrix},
\end{align}
as stated and proven next.

\begin{lemma}
    \label{lem:VdotSpecific_general}
Under the matching condition \eqref{eq:match_MI}, \eqref{eq:condOMEGA}, 
the directional derivative of \eqref{eq:VspecificMI} along \eqref{eq:linearCL} simplifies to
\eqref{eq:Vdot_matchingMI}. Moreover, if  
\eqref{eq:condPosMI1}
% \eqref{eq:LiuAkasaka} 
holds, then $-
\dot V(x)$ in \eqref{eq:Vdot_matchingMI} satisfies item (ii) of Theorem~\ref{theo:definite}.
\end{lemma}

\begin{proof}
First observe that, since $H$ is diagonal, we have
$$
\frac{d}{dt} (\dz^\top H \dz)=\dz^\top H K\dot{x}+\dot{x}^\top K^\top H \dz,
$$
so that we may follow parallel steps to those carried out in the proof of Lemma~\ref{lem:VdotSpecific}, to obtain $-\dot V(x) = 
\frac{1}{2}\smallmat{x \\ \dz}^\top \smallmat{Q_{11} & Q_{12} \\ Q_{12}^\top & Q_{22}}  \smallmat{x \\ \dz}$, with
\begin{align*}
    Q_{11} &= \He ( -(P_0+K^\top H K)(A_0-BK)\\
    &= \He( -P_0A_0+P_0BK -A_0^\top \KT H K ) + \KT \Omega K\\ 
    Q_{12}&= -(P_0+K^\top H K)B+(A_0-BK)^\top K^\top H \\ 
    &= -P_0B + A_0^\top \KT H - \KT \Omega \\ 
    Q_{22} &= \Omega.
\end{align*}
where $\Omega$ is defined in \eqref{eq:condOMEGA}.
Substituting the matching condition \eqref{eq:match_MI} in the identities above leads to the expression of $Q$ given in \eqref{eq:Vdot_matchingMI},
where $S_0$ is defined in \eqref{eq:condPosMI1}.
% Equality \eqref{eq:Q_MIMO_matching} is the multi-input generalization of \eqref{eq:Vdot_matching}.

Consider now the following bound for $Q$, motivated by $S_0\succeq 0$ as per \eqref{eq:condPosMI1},
\begin{align}
\label{eq:QboundLemmaVdotMI}
Q &\succeq     \begin{bmatrix}
  K^\top(2M+\Omega) K  &- K^\top (M+\Omega)\\
-(M+\Omega)K & \Omega
    \end{bmatrix} \\ 
    \nonumber
    &=  \begin{bmatrix}
  K^\top \\
  -I
    \end{bmatrix}
(2M+\Omega)    
     \begin{bmatrix}
  K^\top \\
  -I
    \end{bmatrix}^\top
    +\bigmat{0 &  \KT M \\ M K  & -2M}.
\end{align}
Under the positive definiteness assumption $2M+\Omega\succ 0$
in \eqref{eq:condOMEGA} and due to 
 $M\succeq 0$ diagonal from  \eqref{eq:match_MI}, 
 we obtain that \eqref{eq:QboundLemmaVdotMI}
 implies 
item (ii) of Theorem~\ref{theo:definite} with the selection $R = 2M+\Omega\succ 0$ and $T_0=M \succeq 0$.
%
% As a consequence, all of the assumptions of 
% Theorem~\ref{theo:definite} hold true and 
% the origin is GAS.
\end{proof}

Lemmas~\ref{lem:pos_def_general} and~\ref{lem:VdotSpecific_general} are now exploited in 
% These generalizations are formalized in 
the next theorem, which is our main result for multi-input systems
exploiting the fact that  function $V$ in \eqref{eq:VspecificMI} satisfies the assumptions of Theorem~\ref{th:GAS}.

% \begin{assumption}\label{ass:spectra}The pair $(A_0,K)$ is such that
% $$
% \mathrm{spec}(A_0)\subset \mathbb{C}_{\leq 0},\quad \mathrm{spec}(A_0-BK)\subset \mathbb{C}_{< 0}
% $$
% \end{assumption}
% We can state and prove the following result.

\begin{theo}\label{theo:stabilityMI}
Suppose that symmetric matrices $P_0$ and $H$ satisfy \eqref{eq:condPosMI1}--\eqref{eq:Hdecomposition} and that 
there exists a diagonal matrix $M\succeq0$ satisfying \eqref{eq:match_MI} and \eqref{eq:condOMEGA}.
Then the origin of \eqref{eq:linearCL} is GAS.
%
% . If there exists a decomposition $H=H_1-H_2$ with a diagonal matrix $H_1\succ0$ and a symmetric matrix $H_2\succeq0$ such that
% \begin{equation}\label{eq:H2}
% P_0-K^\top H_2 K\succeq0,
% \end{equation} then function $V$ in \eqref{eq:VspecificMI} is positive definite and radially unbounded.
%
% Suppose also that there exists a diagonal matrix $M\succeq0$ satisfying \eqref{eq:match_MI} and \eqref{eq:condOMEGA}, then the directional derivative of $V$ along \eqref{eq:linearCL} satisfies the second assumption of Theorem~\ref{th:GAS}, namely the origin is GAS.
\end{theo}

\begin{proof}
Combining the results of Lemmas~\ref{lem:pos_def_general} and~\ref{lem:VdotSpecific_general}, we obtain that the conditions of
Theorem~\ref{th:GAS} are satisfied. The result then follows from Theorem~\ref{th:GAS}.
\end{proof}

\newcommand\hsp{\hspace*{1.5cm}}

Theorem~\ref{theo:stabilityMI} establishes analysis conditions for global asymptotic stability of \eqref{eq:linearCL}. An important feature of these conditions is that they are linearly parametrized in the parameters of the Lyapunov certificate \eqref{eq:VspecificMI}. This allows casting the analysis problem into the solution of a linear matrix inequality (LMI). 
To this end, the following formulation allows checking feasibility of the design, while minimizing the condition number of the positive definite matrix $P_0+K^\top H K$.
\begin{align}
\nonumber
&\min_{\substack{\kappa \in \real, P_0 \in {\mathbb S}^m, \Hpos \in {\mathbb D}^m_{>0},\\ \Hneg,M \in {\mathbb D}^m_{\geq 0}}}
\kappa, \text{ subject to }\\
\nonumber
&\hsp  P_0-K^\top \Hneg K \succeq 0, \\
\label{eq:LMI}
&\hsp I \preceq P_0+K^\top (\Hpos - \Hneg) K \preceq \kappa I \\
\nonumber
&\hsp 2M+(\Hpos - \Hneg)KB+B^\top K^\top (\Hpos - \Hneg)\succ 0\\
\nonumber
&\hsp P_0B-A_0^\top K^\top (\Hpos - \Hneg) =K^\top M, \\
\nonumber
&\hsp P_0A_0+A_0^\top P_0\preceq 0.
\end{align}
 Note that restricting 
$P_0+K^\top H K\succeq I$ rather than 
$P_0+K^\top H K\succ 0$ does not introduce any conservatism because all the constraints are homogeneous in the decision variables. 
This LMI-based strategy is summarized in the next corollary.

\begin{coro}
\label{cor:LMI}
Given $A$, $B$ and $K$ of appropriate dimensions,
suppose that the LMI conditions \eqref{eq:LMI}, then the closed-loop \eqref{eq:linearCL} is GAS.
\end{coro}

\section{Case study: downward cart-pendulum}
\label{sec:examples}

\begin{figure}[h!]
\centering
\includegraphics[width=0.7\columnwidth]{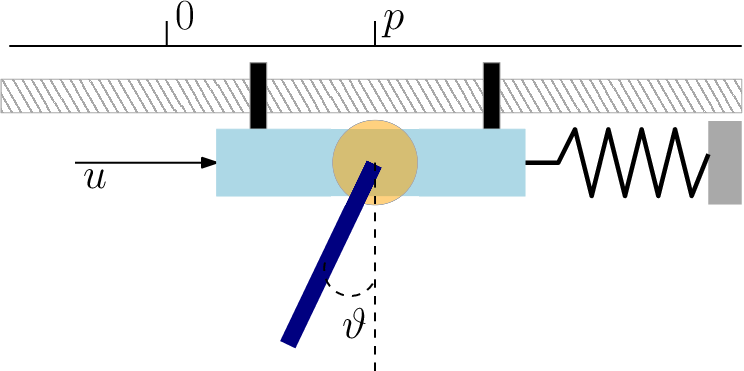}
\caption{Downward spring-cart-pendulum: model.}\label{fig:model}
\end{figure}

We consider the linearized model of the spring-cart-pendulum system taken from \cite{grimm2003antiwindup}, where we assume to control the plant directly through a (saturated) horizontal force input $u$ rather than through a motor as originally done in \cite{grimm2003antiwindup}: 
$$
\begin{array}{rl}
\ddot{p}&=\displaystyle\frac{4}{4M+m}\left(-\kappa p -\frac34 g m\vartheta+\mathrm{sat}(u)\right)\medskip\\
\ddot{\vartheta}&=\displaystyle\frac{3}{4Mr+mr}\left(-\kappa p-g(M+m)\vartheta+\mathrm{sat}(u) \right)
\end{array}
$$
where $\kappa=209.016\, \mathrm{[N/m]}$ is the spring stifness coefficient, $M=0.58\,\mathrm{[Kg]}$ is the mass of the cart, $m=0.21\, \mathrm{[Kg]}$ is the mass of the pendulum, $r=0.305\, \mathrm{[m]}$ is the length of the pendulum and $g=9.81\, \mathrm{[m/s^2]}$ is the gravity acceleration. As illustrated in Figure~\ref{fig:model}, the state $p$ indicates the horizontal displacement of the cart center of mass  
and the state $\vartheta$ indicates the angular position of the pendulum. Defining the full state $x:=\begin{bmatrix} p &\vartheta &\dot{p} &\dot{\vartheta}\end{bmatrix}^\top$,  the plant can be represented by \eqref{eq:plant0} 
with
$$
A_0:=\begin{bmatrix}
0  &   0 &    1 &    0\\
 0   &  0  &   0 &    1\\
  -330.46   &    -2.44&  0   &   0\\
  -812.61 & -30.1   & 0     &    0
\end{bmatrix},\quad B=\begin{bmatrix}
0\\0\\
-1.5810\\
    -3.8878
\end{bmatrix}
$$By direct computation, one can verify that the open-loop state matrix $A_0$ is marginally stable, with two pairs of purely imaginary eigenvalues $\pm\omega_1\mathrm{i}$ and $\pm\omega_2\mathrm{i}$, where $\omega_1=18.3557$ and $\omega_2=4.8641$. 
Aiming at using the developments in Section~\ref{sec:osc} as a guideline, we can consider a coordinate transformation $G\in\mathbb{R}^{4\times 4}$ decoupling the natural modes of the plant, with %\red{I think that there are infinite selections of $G$ giving the expression below. Perhaps we should tell the reader what $G$ we selected and why we did so... }

$$
\tilde{A}_0:=G^{-1}A_0G=\begin{bmatrix}
0&-\omega_1&0 &0\\
\omega_1&0&0&0\\
0&0&0&\omega_2\\
0&0&-\omega_2&0
\end{bmatrix}
$$
and yielding
$$
\tilde{B}:=G^{-1}B=\begin{bmatrix}
b_1\\
0\\
-b_2\\
0
\end{bmatrix},\ b_1=4.1814,\ b_2=0.2936
$$
Setting $\tilde{P}_0:=\mathrm{blkdiag}(b_1^{-1}I_{2\times 2},b_2^{-1}I_{2\times 2})$ and defining a feedback gain $\tilde{K}$ of the form
$$
\tilde{K}:=\displaystyle\begin{bmatrix}k_1& k_2& -\frac{k_1(k_1^2+k_2^2)\omega_1^2}{k_1^2\omega_1^2+k_2^2\omega_2^2}&  \frac{k_2(k_1^2+k_2^2)\omega_1\omega_2}{k_1^2\omega_1^2+k_2^2\omega_2^2}\end{bmatrix}
$$
with $k_1,k_2>0$, the matching condition 
$$
\tilde{P}_0\tilde{B}-h\tilde{A}_0^\top\tilde{K}^\top=\mu\tilde{K}^\top
$$
is satisfied with
\begin{equation}\label{eq:h_example}
h=\frac{k_2}{(k_1^2+k_2^2)\omega_1},\quad \mu=\frac{k_1}{(k_1^2+k_2^2)}
\end{equation}
% \textcolor{teal}{
Resorting to the original coordinates, and defining $K:=\tilde{K}G$, the previous analysis shows that, for any $k_1,k_2>0$, the Lyapunov function 
in \eqref{eq:Vspecific} with  $h$ as in \eqref{eq:h_example} and  
\begin{equation}\label{eq:P0_example}P_0:=G^{-\top} \tilde{P}_0 G^{-1},\quad 
P_0A_0+A_0^\top P_0=0_{4\times 4},
\end{equation} satisfies Lemma~\ref{lem:VdotSpecific}, thereby proving that the saturated closed loop \eqref{eq:linearCL}
is GAS.

\begin{figure}[h!]
\centering
\includegraphics[width=0.9\columnwidth]{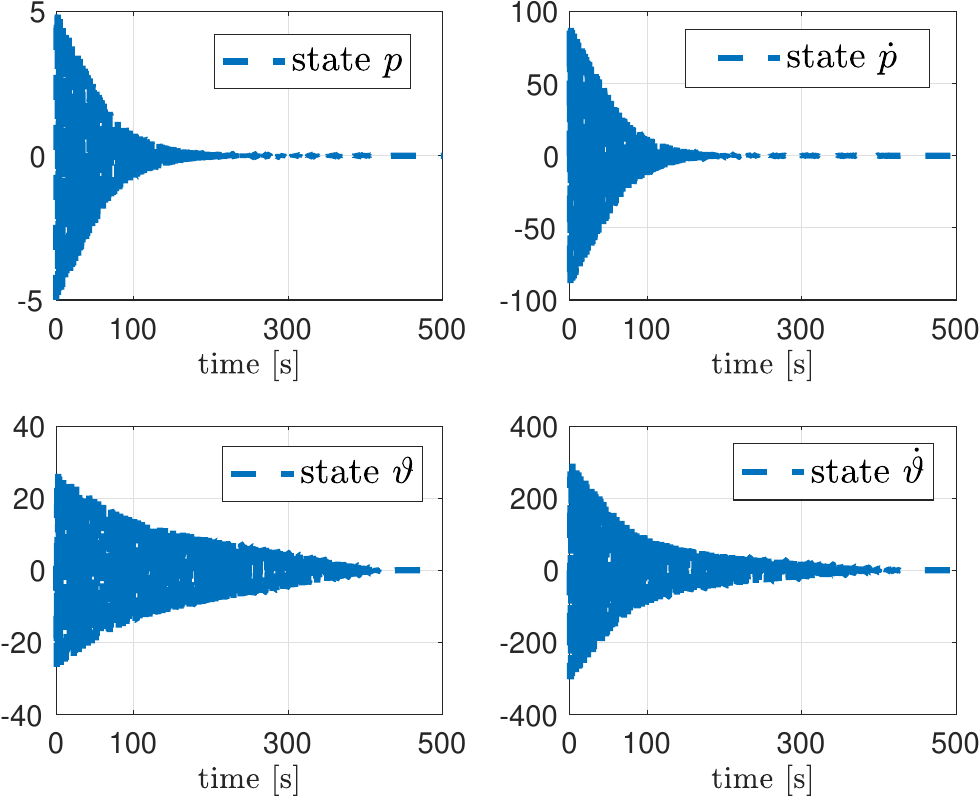}
\caption{Downward cart-pendulum: displacement states $(p,\dot{p})$ (top) and angular states $(\vartheta,\dot\vartheta)$ (bottom).}\label{fig:states}
\end{figure}

\begin{figure}[h!]
\centering
\includegraphics[width=0.9\columnwidth]{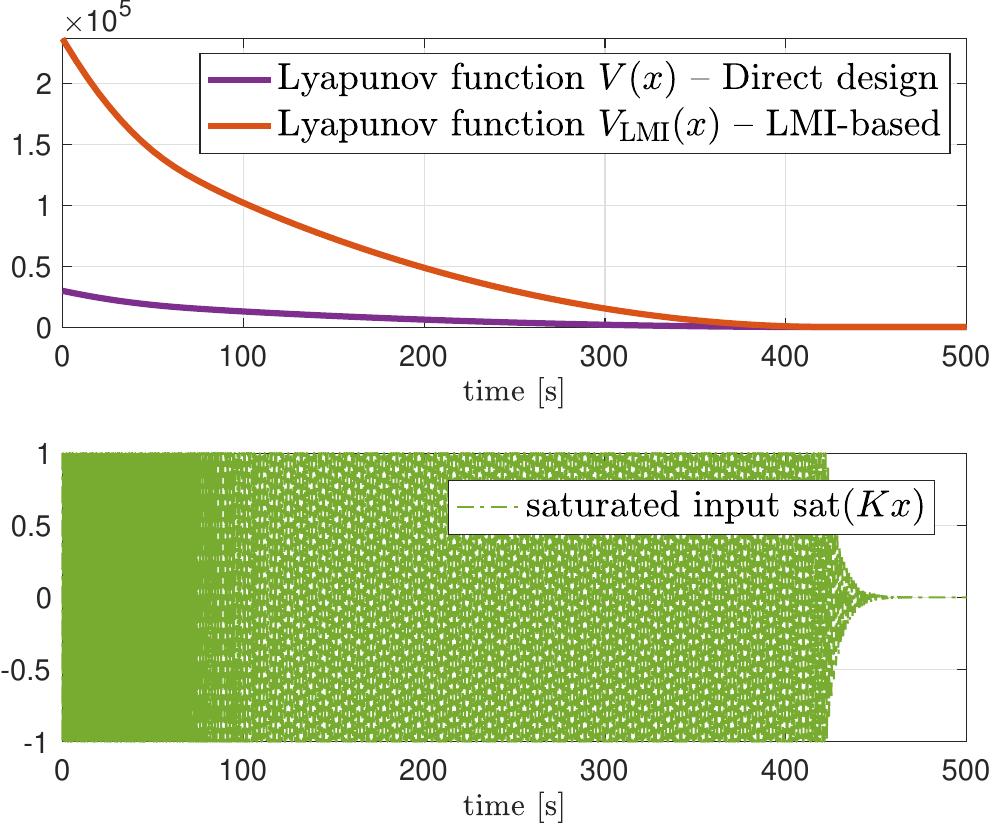}
\caption{Downward cart-pendulum: Lyapunov functions \eqref{eq:Vspecific} 
(top) and saturated input (bottom)}\label{fig:lyap}
\end{figure}

Numerical simulations have been performed with the choice $k_1=1,\ k_2=0.2$. The evolution of the ``displacement'' states 
is depicted in Figure~\ref{fig:states} (top), while the evolution of the ``angular'' states is shown in Figure~\ref{fig:states} (bottom). It can be observed that both are characterized by large oscillations, and that the convergence of the former states is faster. In addition to the Lyapunov function $V$ obtained by direct computation and characterized by the parameters \eqref{eq:h_example}-\eqref{eq:P0_example}, we have considered a second Lyapunov function of the form \eqref{eq:Vspecific}, denoted by $V_{\textrm{LMI}}$, which has been obtained by applying the LMI-based design procedure in Corollary~\ref{cor:LMI}, and characterized by
$$
P_0=10^3\smallmat{
    8.6265 &  -1.5941  & 0 &  0\\
   -1.5941    &0.6151 &  0&    0\\
   0& 0&    0.1951&   -0.0687\\
    0 &  0&-0.0687&    0.0260},\quad h=0.0833 .
$$
Figure~\ref{fig:lyap} (top) reports the evolution of the Lyapunov functions $V(x)$ and $V_{\textrm{LMI}}(x)$, while the bottom plot shows the saturated feedback control $\sat(Kx)$.  In particular, as expected, both functions $V(x)$ and $V_{\textrm{LMI}}(x)$ are monotonically decreasing along the system trajectories, thereby confirming the asymptotic stability. It is also interesting to observe that the saturated control input rapidly oscillates between the upper and the lower limit, thus showing a sort of bang-bang behavior in the initial transient response.

% \begin{figure}[h!]
% \centering
% \includegraphics[width=0.9\columnwidth]{CS-states_theta.eps}
% \caption{Downward cart-pendulum: angular states $(\vartheta,\dot{\vartheta})$}\label{fig:theta}
% \end{figure}

\section{Conclusions}
\label{sec:conclusions}

We leveraged in this paper extended quadratic forms involving the deadzone of the input for obtaining quadratic forms enjoying mild positive definiteness and radial unboundedness properties without quadratic growth in the state. We then use these functions, and the ensuing conditions, to prove global asymptotic stability in cases wehre global exponential stability cannot hold, namely the limit case of ANCBI plants with saturated linear feedback. our construction covers and generalizes Popov-like Lyapunov functions and provides a unified framework for analyzing representative plants, such as the single and double integrator, the integral oscillator and the simple oscillator. We then showed that our construction provides a class of Lyapunov functions leading to convex LMI-based stability conditions for higher order plants involving those non-exponentially converging modes.
Our results have been illustrated on a relevant case study inspired by a mechatronic application.

\appendix 
\section*{Proof of Proposition~\ref{prop:LiuAkasaka}}

% \color{magenta}
To begin with, without loss of generality, we may assume that $A_0$ is in Jordan normal form $J_0$, Indeed, if this is not the case, it suffices to compute a similarity transformation $T$
yielding $T^{-1}A_0T = J_0$, and a matrix $P_0 \succeq 0$ solves \eqref{eq:LiuAkasaka} for $A_0$ if and only matrix $T^\top P_0 T\succeq 0$ solves \eqref{eq:LiuAkasaka} for $J_0$.

To the end of proving necessity, assume that \eqref{eq:LiuAkasaka} holds for some $P_0 \succeq 0$. Denote by ${\mathcal E} \subset \real^n$ the linear subspace of open-loop equilibria,
${\mathcal E} := \{x\in \real^n: A_0x = 0\}$. Denote also by $\Pi =\Pi^{\top}\in \real^{n\times n}$ the orthogonal projection on ${\mathcal E}$, so that the (squared) distance of
a point $x$ from ${\mathcal E}$ corresponds to $|x|^2_{\mathcal E}:=|(I-\Pi) x|^2$.
Using a (e.g., Cholesky) factorization $P_0 = N_0^\top N_0$ we have that
 $x^\top P_0 x = x^\top N_0^\top N_0 x=0$ if and only if $N_0 x=0$,
if and only if $x\in \ker P_0$, which, due to \eqref{eq:LiuAkasaka}
happens only if $A_0x=0$. As a consequence, denoting $V_0:= x^\top P_0 x$,
and exploiting the spectral decomposition of $P_0$, it is immediate to show that there exist
two scalars 
$\underline c$ and $\overline c$ (corresponding to the smallest and the largest nonzero eigenvalues of $P_0$) such that
\begin{align}
\label{eq:sandwich0}
    \underline c |x|^2_{\mathcal E} \leq V_0(x) \leq \overline c |x|^2.
\end{align}
Moreover, due to \eqref{eq:LiuAkasaka},
$\dot V_0(x) = \langle 2P_0x, A_0x\rangle = x^\top (A_0^\top P_0 + P_0 A_0)x \leq 0$
so that $V_0$ is non-increasing along the
the open-loop solutions of $\dot x = A_0 x$. Finally, exploiting \eqref{eq:sandwich0}, we have, for any such solution $x(t)$,
\begin{align}
    |x(t)|^2_{\mathcal E} \leq \tfrac{1}{\underline c}V_0(x(t))\leq \tfrac{1}{\underline c} V_0(x(0))\leq \tfrac{\overline c}{\underline c} |x(0)|^2.
\end{align}
Summarizing, if \eqref{eq:LiuAkasaka} holds for some $P_0 \succeq 0$, then the distance to 
the equilibrium set ${\mathcal E}$ is bounded along all the solutions to $\dot x = A_0 x$. More specifically, with $K_{\mathcal E}>0 := (\tfrac{\overline c}{\underline c})^{\tfrac{1}{2}}$, all solutions satisfy
\begin{align}
\label{eq:xe_bounded}
   |x(t)|_{\mathcal E} \leq K_{\mathcal E}|x(0)|, \quad \forall t\geq 0.
\end{align}

Using \eqref{eq:xe_bounded} we may prove the necessity statement by exploiting the fact that $A_0$ is in Jordan form, which allows grouping distinct eigenvalues in different Jordan blocks, therefore splitting the analysis in a multi-case fashion:\\
(i) First, $A_0$ cannot have eigenvalues with positive real part, because any such eigenvalue generates a non-equilibrium eigenspace associated with exponentially diverging solutions, thereby contradicting \eqref{eq:xe_bounded}.\\
(ii) Secondly, $A_0$ cannot have purely imaginary eigenvalues with multiplicity larger than $1$, because they would generate  a non-equilibrium eigenspace associated with polynomially diverging (oscillatory) solutions, thereby contradicting \eqref{eq:xe_bounded}.\\
(iii) Third, $A_0$ cannot have eigenvalue at the origin with multiplicity larger than $2$, 
such as $A_0 = \begin{smatrix} 0 & I \\ 0 & 0 \end{smatrix}$, with $I$ having dimension larger than $2\times 2$
because in this case we have ${\mathcal E} = \{x\in \real^n: x_1=0\}$, but selecting $x_3(0) = 1$ and all other initial conditions equal to zero, the resulting trajectory provides 
$x_3(t) =1, \;\forall t\geq 0$ and
$x_2(t) = \int_0^t x_3(\tau) d\tau = t, \; \forall t\geq 0$, which diverges polynomially, thereby contradicting \eqref{eq:xe_bounded}.

All of the other possible eigenvalues are allowed in the conditions of Proposition~\ref{prop:LiuAkasaka}, therefore necessity is proven.

Let us now move on to proving sufficiency, and due to the Jordan structure of $A_0$, we can select a block diagonal $P_0$, whose blocks correspond to blocks of $A_0$ where the eigenvalues satisfy the conditions of Proposition~\ref{prop:LiuAkasaka}, as follows:\\
(i) when $A_0$ Hurwitz, it is well known that the Lyapunov inequality $P_0A_0+A_0^\top P_0 \preceq 0$ always admits a positive definite solution $P_0\succ0$.\\
(ii) when $A_0=0$, the trivial selection $P_0=0$ satisfies \eqref{eq:LiuAkasaka}.\\
(iii) when $A_0=\begin{smatrix}0&1\\
0&0   
\end{smatrix}$, constraint \eqref{eq:LiuAkasaka} holds with $P_0=\begin{smatrix}0&0\\
0&1   
\end{smatrix}$ and $S_0=\begin{smatrix}0&0\\
0&0   
\end{smatrix}$, yielding $\ker(P_0)=\ker(A_0)=\mathrm{span}(\begin{smatrix}1\\
0   
\end{smatrix})$. \\
(iv) when $A_0=J_\omega:=\begin{smatrix}0&\omega\\
-\omega&0   
\end{smatrix}$ with $\omega\in\mathbb{R}$, constraint \eqref{eq:LiuAkasaka} 
holds with the selections $P_0=\begin{smatrix}1&0\\
0&1   
\end{smatrix}$ and $S_0=\begin{smatrix}0&0\\
0&0   
\end{smatrix}$, where $\ker(P_0)=\ker(A_0)=\{0\}$.

Since all of the possible eigenvalue selections have been covered, this completes the sufficiency of the proof.

\bibliographystyle{plain}

\bibliography{biblio}

\end{document}